\documentstyle[12pt]{article}
\title{Shape curvatures  and  transversal fluctuations  in the first passage percolation model
} 
\author{Yu Zhang\footnote{Research supported by NSF grant DMS-0405150. }}
\date{}

\begin{document}
\baselineskip .20in
\maketitle
\begin{abstract}
We consider the first passage percolation model on ${\bf Z}^2$.
In this model,  
$\{t(e): e\mbox{ an edge of }{\bf Z}^2 \}$ is an independent identically distributed family with 
a common distribution $F$.
We denote by
$T({\bf 0}, v)$  the passage time from the origin to $v$ for $v\in {\bf R}^2$ and
$B(t)=\{v\in {\bf R}^d: T({\bf 0}, v)\leq t\}.$
It is well known that if $F(0) < p_c$, there exists a compact shape ${\bf B}_F\subset {\bf R}^2$ such that
for all $\epsilon >0$,
$t {\bf B}_F(1-\epsilon )  \subset {B(t) } \subset t{\bf B}_F(1+\epsilon )$, eventually with a probability 1.
For each shape boundary  point $u$, we denote its  right- and left-  curvature exponents
by $\kappa^+(u)$ and $\kappa^-(u)$.
In addition, for each  vector $u$, we  denote  the transversal fluctuation  exponent by $\xi(u)$.
In this paper,  we can show that $\xi(u) \leq 1-\max\{\kappa^-(u)/2, \kappa^+(u)/2\}$ for all shape boundary
points $u$.

To pursue a curvature on ${\bf B}_F$, we consider passage times with
 a special distribution infsupp$(F)=l$ and $F(l)=p > \vec{p}_c$, where $l$ is a positive number  and 
$\vec{p}_c$ is a critical point for the oriented percolation  model.  With this distribution,
it is known that there is a flat segment on the shape boundary between angles 
$0< \theta_p^- < \theta_p^+< 90^\circ$. 
In this paper, we show that  the shape are strictly convex at the directions $\theta_p^\pm$.
Moreover, we also show that for all $r>0$,
$\xi((r, \theta^\pm_p)) = 0.5$ and 
$\xi((r, \theta)) =1$ for all  $\theta_p^- <\theta< \theta_p^+$ and $r>0$.
Note that this rules out the conjecture  that $\xi(u)=2/3$ for all $u$.
In addition, if  $\chi(u)$ is the longitudinal exponent, it is believed that $\chi(u)= 2\xi(u)-1$. Furthermore,
 it is estimated that
$\chi(u) \geq (1-\xi(u))/2$ for $F(0) < p_c$ and infsupp$(F)=0$. 
However, both the equation and the inequality do not hold
for our special distribution when $\theta=\theta^\pm_p$.

\end{abstract}
\section {Introduction of the model and results.}
The first passage percolation model was introduced  in 1965 by Hammersley and  Welsh. 
In this model, we consider the ${\bf Z}^2$ lattice as a graph with edges connecting each pair of vertices
$u_1=(r_1, \theta_1)$ and $u_2=(r_2,\theta_2)$ 
with $\|u_1- u_2\|=1$, where $\|u_1-u_2\|$ is the Euclidean distance between $u_1$ and $u_2$.
In this paper, we always use the polar coordinates $\{(r, \theta)\}$, where $r$ and $\theta$ present the radius 
and the  angle between the radius and the $X$-axis, respectively.

We assign independently to each edge a non-negative {\em passage time} $t(e)$ with a common distribution $F$.
More formally, we consider the following probability space. As the
sample space, we take $\Omega=\prod_{e\in {\bf Z}^2} [0,\infty),$ whose points 
are called {\em configurations}.
Let ${\bf P}=\prod_{e\in {\bf Z}^2} \mu_e$ be the corresponding product measure on $\Omega$, where
$\mu_e$ is the measure on $[0, \infty)$ with distribution $F$. The
expectation  with respect to $P$ is denoted by ${\bf E}(\cdot)$.
For any two vertices $u$ and $v$, 
a path $\gamma$ from $u$ to $v$ is an alternating sequence 
$(v_0, e_1, v_1,...,v_i, e_{i+1}, v_{i+1},...,v_{n-1},e_n, v_n)$ of vertices $v_i$ and 
 edges $e_i$ between $v_i$ and $v_{i+1}$ in ${\bf Z}^2$ with $v_0=u$ and $ v_n=v$. 
Given such a path $\gamma$, we define its passage time  as 
$$T(\gamma )= \sum_{i=1}^{n} t(e_i).\eqno{}$$
For any two sets $A$ and $B$, we define the passage time from $A$ to $B$
as
$$T(A,B)=\inf \{ T(\gamma)\},$$
where the infimum is over all possible finite paths from some vertex in $A$ to some vertex in $B$.
A path $\gamma$ from $A$ to $B$ with $T(\gamma)=T(A, B)$ is called the {\em optimal path} of 
$T(A, B)$. The existence of such an optimal path has been proven (see Kesten (1986)).
We also want to point out that the optimal path may  not be  unique.
If we focus on a special configuration $\omega$, we may write $T(A, B)(\omega)$
instead of $T(A, B)$. 
When $A=\{u\}$ and $B=\{v\}$ are single vertex sets, $T(u,v)$ is the passage time
from $u$ to $v$. We may extend the passage time over ${\bf R}^2$.
If $x$ and $y$ are in ${\bf R}^2$, we define $T(x, y)=T(x', y')$, where
$x'$ (resp., $y'$) is the nearest neighbor of $x$ (resp., $y$) in
${\bf Z}^2$. Possible indetermination can be eliminated by choosing an
order on the vertices of ${\bf Z}^2$ and taking the smallest nearest
neighbor for this order.

With these definitions, we would like to introduce the basic developments and questions
in this field.  Hammersley and  Welsh (1965) first studied the  point-point  and the point-line passage times
defined as follows:
$$a_{m,n}= \inf\{T(\gamma ): \gamma  \mbox{ is a path from $(m, 0)$ to  $(n, 0)$}\},
$$
$$b_{m,n}= \inf\{T(\gamma ): \gamma  \mbox{ is a path from $(m, 0)$ to  $\{x=n\}$}\}.
$$
It is well known (see Smythe and Wierman (1978)) that if $Et(e) < \infty$,
$$
\lim_{n\rightarrow \infty}{1\over n} a_{0,n}= \lim_{n\rightarrow \infty}{1\over n} b_{0,n}= \mu \mbox{ a.s. and in } L_1,\eqno{(1.1)}
$$
where the non-random constant $\mu=\mu(F)$ is  called the {\em time constant}. 
Later, Kesten showed (see Theorem 6.1 in Kesten (1986)) that
$$\mu=0\mbox{ iff } F(0)\geq p_c,\eqno{(1.2)}$$
where $p_c=1/2$ is the critical probability for Bernoulli (bond) percolation on ${\bf Z}^2$.

Given a non-zero vector $(r, \theta)\in {\bf R}^2$, by the same arguments as in (1.2) and (1.3), if $Et(e) < \infty$, then
$$
\lim_{n\rightarrow \infty}{1\over n} T({\bf 0}, (nr, \theta)) = \lim_{n\rightarrow \infty} {1\over n} E T({\bf 0}, (nr, \theta))=\mu((r, \theta)) \mbox{ a.s. and in } L_1,\eqno{(1.3)}
$$
and 
$$\mu((r, \theta))=0\mbox{ iff }F(0) \geq p_c.$$
For  convenience, we assume that $t(e)$ is not a constant and satisfies strong moment requirement in this paper:
$$\int e^{\lambda x} dF(x) < \infty \mbox{ for some } \lambda >0.\eqno{(1.4)}$$

When $F(0) < p_c$, the map $x \rightarrow \mu(x)$ induces  a norm on
${\bf R}^2$. The unit radius ball for this norm is denoted by ${\bf B}:={\bf B}(F)$
and is called the {\em asymptotic shape}. The boundary of ${\bf B}$ is
$$\partial {\bf B}:= \{ x \in {\bf R}^2: \mu(x)=1\}.$$
${\bf B}$ is a compact convex
deterministic set and  $\partial {\bf B}$  is a continuous convex closed
curve (Kesten (1986)).  Define for all $t> 0$,
$$B(t):= \{v\in {\bf R}^2, \ T( {\bf 0}, v) \leq t\}.$$
The shape theorem (see Theorem 1.7 of Kesten (1986)) is the well-known result stating that for any
$\epsilon >0$,
$$t{\bf B}(1-\epsilon)  \subset {B(t) } \subset t{\bf B}(1+\epsilon )
\mbox{ eventually w.p.1.}\eqno{}$$

The  most difficult problem aspect in this field  is to question the {\em transversal fluctuations} of optimal paths (see Hammersley and Welsh (1965), Kingman (1973), Smythe and Wierman (1978), Kesten (1986), and Newman and Piza (1995)). 
Let us use the notations  of Newman and Piza (1995) to define transversal fluctuations. 
Given a vector $u=(r, \theta)$, let $M_n({u})$ denote the random set
of all vertices in ${\bf Z}^2$ belonging to some optimal paths of $T({\bf 0}, n u)$. 
Let ${\bf L}_\theta$ denote the line passing the origin with the angle $\theta$ between the line and  the $X$-axis.
Let the  transversal fluctuations of optimal paths be
denoted by
$$h_n({u})=\max\left\{{\bf dist}(v, {\bf L}_\theta) : v\in M_n({u})\right\},$$
where
$$\mbox{{\bf dist}}(A, B)=\min\left\{\|u-v\|: u\in A, v\in B\right\}$$
for some sets $A, B\subset {\bf R}^2$.
With this definition, when $u=(1,0)$, Hammersley and Welsh (1965) asked:\\
$$\mbox{ does } h_n({u})/n \rightarrow 0 \mbox{?}$$
This well-established conjecture is called the {\em height conjecture}.
Let  the {\em transversal fluctuation exponent} $\xi({u})$ be  defined as
$$\xi(u)=\inf \left\{ s>0:  {\bf P}[h_n({ u} )\leq n^s ]\geq C \mbox{ for all large $n$}\right\}.$$
In this paper,  $C$ and $C_i$ are always positive  constants that may depend on $F$, $\delta$, $p$, or other parameters,
but not on $t$, $m$, or $n$. 
Their values 
are not significant and  change from 
appearance to appearance.
It is well-known  (see Theorem 8.17 in Smythe and Wierman (1978)) that for all $u$,
$$\xi(u) \leq 1.$$
The main conjecture is to ask, for some $u$, 
$$\xi({u}) < 1.\eqno{(1.5)}$$

This problem has also received a great amount of attention  from 
 statistical physicists because of its equivalence with one version of the Eden growth model (see  Krug and Spohn (1992)).
Statistical physicists believe that   
$$\xi(u)=2/3 \mbox{ for all }u.\eqno{(1.6)}$$
They also believe that $\xi(u)$ will   decrease  as dimensions increase.
For some growth models, as the increasing subsequence model (see Johansson (2000)), $\xi(u)$ is indeed
$2/3$ for $u=(r, \pi/4)$.

Mathematicians have also made significant efforts in this direction. 
Perhaps Kesten (1986) first noticed that the transversal fluctuations might depend on the behavior of the curve of
$ \partial {\bf B}_F$ around $u$ (see his book, pages 262 - 263 (1986)).
Later, Newman and Piza (1995)  explored a deeper result in this direction. They found out that the quantity
$\xi(u) $ depends on the curvature of ${\bf B}_F$ at the direction of $u$. More precisely, as they defined,
$u\in \partial {\bf B}_F$ is said to be strictly convex if the following occurs:

1. There is a subset $S'$ of the circuit boundary of  $\partial D$ that is open and that contains $u$,
but does not contain a point of ${\bf B}_F^\circ$.

2. $D\cap {\bf B}^\circ_F\neq \emptyset$,
where $A^\circ$ is the interior of $A$ for some set $A$ in ${\bf R}^2$.

With this definition, Newman and Piza (1995) showed in their Theorem 6 that if $F(0) < p_c$ and $ \partial {\bf B}_F$ is 
strictly convex at $u$, then
 $$\xi(u)\leq 3/4.\eqno{(1.7)}$$
Since ${\bf B}_F$ is convex, $\partial {\bf B}_F$ has a strictly convex point. 
\begin{figure}\label{F:graph}
\begin{center}
\setlength{\unitlength}{0.0125in}%
\begin{picture}(200,330)(67,850)
\thicklines
\put(0,1100){\oval(60,60)}
\put(14,1129){\circle*{4}}
\put(30,1115){\circle*{4}}
\put(35,1107){$r_{{\bf B}}(\theta)$}

\put(0, 1100){\line(1,2){50}}
\put(0, 1100){\line(2,1){80}}

\put(48,1124){\circle*{4}}

\put(-10,1133){\line(6,-1){75}}
\put(70,1114){${{{\bf S}_{\theta_0}^+}}$}
\put(45,1205){${{{\theta_0}}}$}
\put(80,1135){${{{\theta}}}$}
\put(30,1130){${d^+_\theta}$}

\put(-20,1080){${\partial {\bf B}_F}$}
\put(-10,1050){$\mbox{ (a)}$}

\put(300,1100){\oval(60,60)}

\put(315,1131){\circle*{4}}
\put(260,1131){$r_{{\bf B}}(\theta)$}

\put(293,1130){\circle*{4}}

\put(293,1142){\circle*{4}}

\put(300, 1100){\line(1,2){50}}
\put(280, 1220){\line(1,-6){20}}

\put(270,1153){\line(2,-1){75}}
\put(350,1114){${{{\bf S}_{\theta_0}^-}}$}
\put(345,1205){${{{\theta_0}}}$}
\put(280,1220){${{{\theta}}}$}
\put(300,1140){$d_\theta^-$}

\put(280,1080){${\partial {\bf B}_F}$}
\put(290,1050){$\mbox{ (b)}$}
\end{picture}
\end{center}
\vskip -2in
\caption{{\em The graph in (a) shows that $\partial {\bf B}_L$ has a right $l$-curvature. 
The distance from $r_{\bf B}(\theta)$ to ${\bf S}^+_{\theta_0}$ is larger than  $ C(d^+_\theta)^{1/l}$. 
The graph in (b) shows that $\partial {\bf B}_L$ has a left $l$-curvature. 
The distance from $r_{\bf B}(\theta)$ to ${\bf S}^-_{\theta_0}$ is larger than  $ C(d^-_\theta)^{1/l}$. }}
\end{figure}
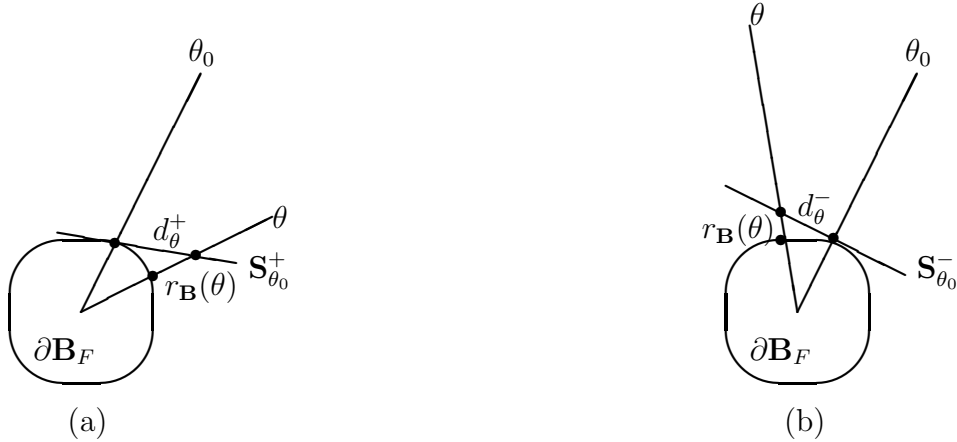

However, the curve on $\partial {\bf B}_F$ can be much more complicated. For example, it is known that
there is a flat segment on the curve of $\partial {\bf B}_F$ for some $F$. Some points in
$\partial {\bf B}_F$ are neither strictly convex nor  flat (see the following examples). 
In these situations, we need  a more general definition for curvatures. Furthermore, 
we might also ask what the behavior of $\xi(u)$ is for this more general definition.  
Intuitively, as a standard definition for  curvature,
the curve $f(x)=|x|^m$  is strictly convex at $x=0$ only if $0<m\leq 2$.
However, to give a more general and quantitative definition, we may say that $f(x)$ has an $l=1/m$-curvature at $x=0$. 
If $f(x)$ is asymmetric at $u$, we may also consider its curvature from  the left and right at $u$.
For example,  the curve
$$f(x)=\left\{\begin{array}{ll}
x^m &\mbox{ if $x\geq 0,$}\\
0 &\mbox{ if $ x< 0,$}
\end{array} \right.$$
has a right $l=1/m$-curvature and a left $0$-curvature at $x=0$.
With this motivation in mind, 
let us  present a precise  definition for the $l$-curvature  at the boundary of the shape. 
We denote by $ \partial {\bf B}_F=\{(r_{\bf B}(\theta), \theta)\}$.
For $(r_0, \theta_0)\in \partial {\bf B}_F$ 
with $0\leq \theta_0\leq \pi/2$, since ${\bf B}_F$ is convex, there exists a line passing through
$\partial {\bf B}_F$ at $(r_0, \theta_0)$, and ${\bf B}_F$ is below the line (see Fig. 1). 
We denote the line by
$${\bf S}_{\theta_0} =\{(r_{\bf S}(\theta),\theta)\}\eqno{(1.8)}$$
for $r_{\bf B}(\theta) \leq r_{\bf S}(\theta)$ and $r_{\bf B}(\theta_0) = r_{\bf S}(\theta_0)$.

Note that there might be many such lines $\{{\bf S}_{\theta_0}\}$ passing through $(r_0, \theta_0)$.
For each line, there is an angle between the line and  the $X$-axis. 
We select two lines with the largest and the smallest angles  from these lines, denoted by
${\bf S}^+_{\theta_0}$ and ${\bf S}^-_{\theta_0}$ (see Fig. 1).  
If ${\bf S}_{\theta_0}$ is unique, then
${\bf S}_{\theta_0}^+={\bf S}_{\theta_0}^-$. 
We  denote the two lines with polar coordinates by
$${\bf S}^+_{\theta_0}=\{(r_{\bf S}^+(\theta), \theta)\}\mbox{ and } {\bf S}^-_{\theta_0}=\{(r_{\bf S}^-(\theta), \theta)\}.$$
With points on $(r_{\bf S}^+(\theta), \theta)\in {\bf S}^+_{\theta_0}$ or  $(r_{\bf S}^-(\theta), \theta)\in {\bf S}^-_{\theta_0}$, let 
$$d_\theta^+=\|(r_{\bf S}^+(\theta), \theta)-(r_0, \theta_0)\|\mbox{ and } d_\theta^-=\|(r_{\bf S}^-(\theta), \theta)- (r_0, \theta_0)\|,$$
the distances from the left and the right points on the line to $(r_0, \theta_0)$ (see Fig. 1).
Clearly,
$$d_{\theta_0}^-=d^+_{\theta_0}=0.$$
On the other hand, by the continuity of $\partial {\bf B}_F$,
we always select $\delta$ such that $d^\pm_\theta < 1.$

With these definitions, we say $(r_0,\theta_0)\in \partial {\bf B}_F$ has at least a left  $l$-curvature 
for $0\leq l \leq 1$ if,
for $\theta\in [\theta_0-\delta, \theta_0]$ with some small $\delta >0$,
$$r_{\bf S}^-(\theta)- r_{\bf B}(\theta) \geq C(d^-_\theta)^{1/l}\eqno{(1.9)} $$
for some constant $C$ that does not depend on $\theta$, where $(r_{\bf B}(\theta),\theta)\in \partial {\bf B}_F$
and $(r_{\bf S}^-(\theta), \theta)\in {\bf S}^-_{\theta_0}$ defined above.
Since $\partial {\bf B}_F$ is convex, any $u\in \partial {\bf B}_F$ has at least a zero curvature.
We denote by
$$\kappa^-(u)=\sup\{l: u \mbox{ has at least a left $l$-curvature}\}$$
the left curvature of $\partial {\bf B}_F$ at $u$.
Similarly, we say $(r_0,\theta_0)\in {\bf B}$ has at least a right $l$-curvature if, for $\theta\in [\theta_0, \theta_0+\delta]$,
$$ r_{\bf S}^+(\theta)- r_{\bf B}(\theta)\geq C(d^+_\theta)^{1/l}.\eqno{(1.10)}$$
We denote by
$$\kappa^+(u)=\sup\{l: u \mbox{ has at least a right $l$-curvature}\}$$
the right curvature of $\partial {\bf B}_F$ at $u$.

Similarly, we can define the curvatures for all $\theta$.
With these definitions, we have a few remarks.\\

{\bf Remark 1.} By our definition, $\kappa^\pm(u) \leq 1$ for  all $u\in \partial {\bf B}_F$. We say $u$ is a {\em sharp point}
if $\kappa^\pm(u)=1$. We also say $u$ is  a {\em right-} or {\em left-flat} if $\kappa^\pm(u)=0$.\\  

{\bf Remark 2.} We may 
replace the inequalities in (1.9) and (1.10)  by 
$$r_{\bf S}^\pm(\theta)- r_{\bf B}(\theta)\geq g(\theta) (d^\pm_\theta)^{1/l}$$
for a slow growth function $g$. However, for simplicity, we will not attempt  this method.\\

{\bf Remark 3.}
It can be verified that if a point of $\partial {\bf B}_F$ is  strictly convex, defined by Newman and Piza (1995)), 
then $\kappa^\pm(u)\geq 0.5$. We say $u$ is strictly {\em right-} or {\em  left-convex} at $u$ if 
$\kappa^\pm \geq 0.5$.\\ 

With this weaker version of curvature, we show  the following theorem.\\

{\bf Theorem 1.} {\em If $F$ satisfies $F(0) < p_c$ and the tail assumption in (1.4), then for 
$u\in \partial {\bf B}_F$ with $q=\max\{\kappa^-(u), \kappa^+(u)\}$, and for $\delta >0$,
there exist constants $C_i=C_i(F, \delta)$ for $i=1,2,3$ such that}
$${\bf P}\left [ { h_n(u) \over  n^{1-q/2+\delta}} > C_1\right] \leq C_2 \exp(-C_2 n^{\delta/2}).\eqno{(1.11)}$$

With Theorem 1, the following corollary is a direct consequence of Theorem 1 and the definition of $\xi(u)$.\\

{\bf Corollary.} {\em Under the same assumption of Theorem 1,}
$$\xi(u) \leq 1-\max\{\kappa^-(u)/2, \kappa^+(u)/2\}.\eqno{(1.12)}$$

{\bf Remark 4.} The proof of Theorem 1 depends on the following estimates (Alexander (1997) and Kesten (1993)).
For $\delta >0$, there exist $C_i=C_i(F, \delta)$ for $i=1,2$ such that 
$${\bf P}\left[|T({\bf 0}, nu) -n\mu_F(u)|\geq n^{1/2+\delta}\right]\leq C_1\exp(-C_2 n^\delta).$$
It is believed that
$${\bf P}\left[|T({\bf 0}, nu) -n\mu_F(u)|\geq n^{1/3+\delta}\right]\leq C_1\exp(-C_2 n^\delta).\eqno{(1.13)}$$
If (1.13) indeed holds, then we can show a better result in Theorem 1 to have
$$\xi(u) \leq 1-\max\{2\kappa^-(u)/3, 2\kappa^+(u)/3\}.$$

{\bf Remark 5.} When $u=(r, 0)$, by  symmetry, $\kappa^-(u)=\kappa^+(u)$. We conjectured that
$$\kappa^-((r, 0))=\kappa^+((r,0))=1/2\mbox{ and } \xi(u)=2/3.$$

With Theorem 1, we may ask whether $0<\kappa^\pm (u) $ for  $u\in \partial {\bf B}_F$.
If it is, what is $\kappa^\pm (u)$? Indeed, one of the major conjectures for the shape is to ask
(see Howard's conjecture 6 (2000)): can any particular direction be shown  to 
be strictly convex for any non-trivial $F$?   To investigate the question, 
we would like to introduce a few results in the oriented percolation 
model. Recall that the classical grid ${\cal L}$ for oriented percolation
is given by ${\cal L}:= \{(m,n) \in {\bf Z} ^2:  m + n  \mbox{ has even parity}, n \geq 0 \}$. 
 Thus, ${\cal L}$ is ${\bf Z}^2$ rotated
by $\pi/4$ and correctly dilated. The edges in ${\cal L}$ are
from $(m,n) \in {\cal L}$ to $(m+1, n+1)$ and to $(m - 1,
n + 1)$.  To each edge $e$ we assign a passage
time 1 with probability $F(1)=p$ and a time larger than 1 with
probability $1 - p$.

We consider all paths starting from
 $\{(x,y) \subset {\bf Z}^2: x \leq 0, \ y = 0 \}$
in the {\em oriented} graph using  $n$ type 1 oriented edges in
${\cal L}$ and let
 $(r_n(p),n)$ denote the right-most point (right-hand edge) of all such paths.
We will often simply refer to the scalar $r_n(p)$ as the
right-hand edge.
By a subadditive argument (see Durrett (1984)), it is well known that the right-most point $(r_n(p),n)$  satisfies
\begin{eqnarray*}
\lim_{n\rightarrow \infty} \frac{r_n(p)}{n} = \alpha_p \ \ \
\mbox{ a.s. and  in  }L^1,
\end{eqnarray*}
Let critical probability be
$$\vec{p}_c=\inf\{p: \alpha_p >0\}.$$
It is also well known (see Durrett (1984)) that $0< \vec{p}_c < 1.$ 
 When $p > \vec{p}_c$, $\alpha_p\in (0,1]$
 is  called the asymptotic speed of
super-critical oriented percolation on the edges of ${\cal L}$. It
describes the drift of the right-most point at level $n$.

With this model, we focus on  the following special distributions
investigated by Durrett and Liggett (1981):
$$\mbox{infsupp}(F)=l>0 \mbox{ and } F(l)=p > \vec{p}_c. \eqno{}$$
Clearly, $F(0)=0 < p_c$, and the shape ${\bf B}_F$ is compact. Without loss of generality,
we may replace $l$ by 1 to have
$$\mbox{infsupp}(F)=1>0 \mbox{ and } F(1)=p > \vec{p}_c. \eqno{(1.14)}$$

With this special distribution,
Durrett and Liggett (1981) found that  shape ${\bf B}_F$ contains a flat segment on the diagonal direction.
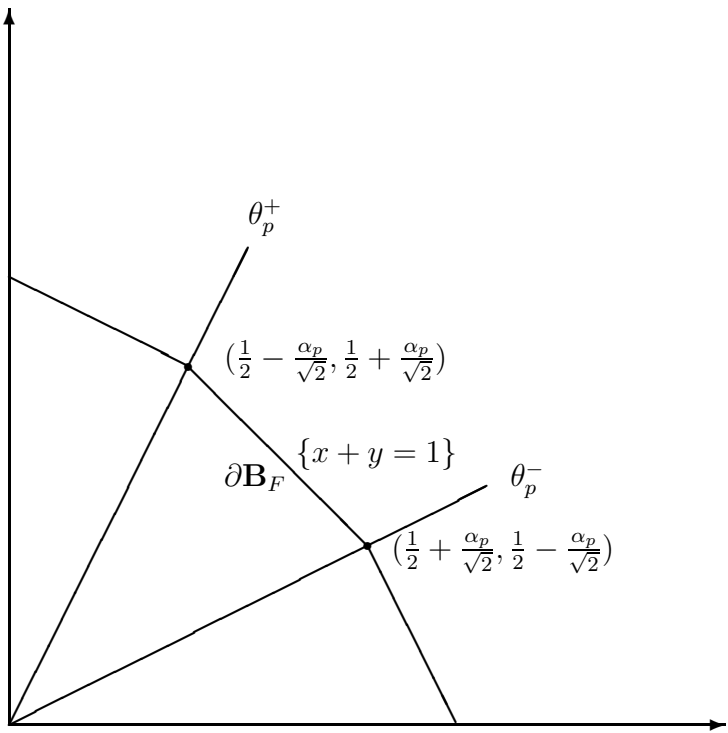
\begin{figure}\label{F:graph}
\begin{center}
\setlength{\unitlength}{0.0125in}%
\begin{picture}(200,230)(67,800)
\thicklines
\put(0, 800){\vector(1,0){300}}
\put(0, 800){\vector(0,1){300}}
\put(0,988){\line(2,-1){75}}
\put(75,950){\circle*{3}}
\put(150,875){\circle*{3}}
\put(0, 800){\line(1,2){100}}
\put(0, 800){\line(2,1){200}}
\put(75,950){\line(1,-1){75}}
\put(150,875){\line(1,-2){37}}
\put(90, 950){${({1\over 2}-{\alpha_p\over \sqrt{2}}, {1\over 2}+{\alpha_p\over \sqrt{2}})}$}
\put(160, 870){${({1\over 2}+{\alpha_p\over \sqrt{2}}, {1\over 2}-{\alpha_p\over \sqrt{2}})}$}
\put(210,900){${{\theta_p^-}}$}
\put(100,1010){$\theta_p^+$}
\put(120,910){$\{x+y=1\}$}
\put(90,900){$\partial {\bf B}_F$}

\end{picture}
\end{center}
\caption{{\em The graph shows that shape ${\bf B}_F$ contains  a flat segment.}}
\end{figure}

Later, Marchand (2002) presented the  precise locations of
the flat segment in the shape.
More precisely,  two polar coordinates in the first quadrant are denoted by $(\sqrt{1/2+ \alpha^2_p}, \theta_p^\pm)$
(see Fig. 2), where
$$\theta_p^{\pm} =\arctan \left( { 1/2\mp \alpha_p/\sqrt{2} \over  1/2\pm \alpha_p/\sqrt{2} }\right ),$$
and $\alpha_p\geq 0$ is the asymptotic speed defined above. 
Note that $\theta_p^- < \theta_p^+$, if $p > \vec{p}_c$, and that $\theta_p^-=\theta_p^+$, if $p=\vec{p_c}$. 
 Marchand (see Theorem 1.3 of Marchand (2002)) showed that under (1.14),
$$\partial {\bf B}_F\cap \{(x,y)\in {\bf R}^2, |x|+|y|=1\}=\mbox{ the segment from $(\sqrt{1/2-\alpha^2_p}, \theta_p^+)$ to $ (\sqrt{1/2+\alpha^2_p},\theta_p^-)$},\eqno{}$$
where the segment will shrink as a point $(1/\sqrt{2}, \pi/4)$ when $p=\vec{p}_c$.
This segment is called the {\em flat segment} of shape ${\bf B}_F$. The cone between $\theta_p^-$ and $\theta_p^+$ 
is called the {\em oriented percolation cone}. 

To understand why this is called the oriented percolation cone, we introduce the 
following oriented path. Let us define the northeast-oriented paths. 
A path (not necessarily  a $1$-path) is said to be a  northeast
  path, called an NE path, if each vertex $u$  of the path  has only one existing edge,
either from $u$ to $u+(1,0)$ or  to $u+(0,1)$. 
If there exists an NE path from $u$ to $v$, then there exists a southwest path, called an SW path, from $v$ to $u$.
We denote by $u\rightarrow v$ if there exists an NE $1$-path
from $u$ to $v$.
 In fact, if we rotate ${\bf Z}^2$ $45^\circ$, then northeast paths will be the oriented paths defined on ${\cal L}.$
In particular, if we want to emphasize an NE $1$-path $\gamma$ from $u$ to $v$, we write
$u \stackrel{\gamma}{\rightarrow} v$. We also denote by $A\rightarrow B$ for two sets $A$ and $B$ if there exist
$u\in A$ and $v\in B$ such that $u\rightarrow v$.

For any vector $(r, \theta)$ with $\theta_1 \leq \theta \leq \theta_2$, under (1.14),  with a positive probability
there is a northeast path $\gamma$ from $(0,0)$ to $(nr, \theta)$ with only 1-edges
 (see (3.2) in Yukich and Zhang (2006)), so we call the cone  between $\theta_1$ and $\theta_2$ 
the oriented percolation cone. 
With these definitions, we may investigate  the transversal fluctuations in an oriented percolation cone.\\

{\bf Theorem 2.} {\em For a  vector $u=(r, \theta)$ with $\theta_p^- < \theta < \theta_p^+$, if $F$ satisfies (1.14), then}
 $$\xi(u)=1.\eqno{(1.15)}$$

By Theorem 2, the conjecture in (1.6) cannot be true. However, the following theorem is more surprising since
 it shows that $h_n(u)$ is normally diffusive for some $u$. \\

{\bf Theorem 3.} {\em For a  vector $u=(r, \theta^\pm_p)$, if $F$ satisfies (1.14), then}
$$\xi(u)=0.5.\eqno{(1.16)}$$

{\bf Remark 6.} When $u=(r, \theta^\pm_p)$, it follows from Theorem 3 that
the transversal fluctuation is normally diffusive. In fact, we can show the following   stronger estimates.
There exist $\delta=\delta(F, u) >0$ and $C_i=C_i(F, \delta)$ for $i=1,2$ such that
$${\bf P}\left [{h_n(u)} \leq \delta n  \right]\leq C_1\exp(-C_2 n), \eqno{(1.17)}$$
when $u=(r, \theta)$ for $\theta_p^- < \theta < \theta_p^+$, and
$${\bf P} \left [{h_n(u)} \geq n^{1/2+\delta} \right ] \leq C_1\exp(-C_2 n^\delta), \eqno{(1.18)}$$
when $u=(r, \theta^\pm_p)$.
For $u=(r, \theta)$ with $\theta < \theta_p^-$ or $\theta > \theta_p^+$, we  conjecture that
$$\xi(u) =2/3.$$

As we mentioned before, another interesting and  important question is how to find a curved point on 
$\partial {\bf B}_F$.
In this paper, we find that  there is a particular curved point on $\partial {\bf B}_F$.\\

{\bf Theorem 4.} {\em If $F$ satisfies (1.4) and (1.14), then} 
$$\kappa^+(r, \theta^-_p)\geq 0.5\mbox{ and } \kappa^-(r, \theta^+_p)\geq 0.5,$$
for $(r, \theta^\pm)\in \partial {\bf B}_F$.\\

{\bf Remark 7.} Theorem 4 tells us that $(r, \theta^\pm_p)\in \partial {\bf B}_F$ is 
strictly left or right-convex., so we partially answer the conjecture asked by Howard (2000).
We believe that $\kappa^\pm(u)=0$, when $u=(r, \theta_p^\pm)\in \partial {\bf B}_F$ with  $F$ 
satisfying (1.4) and (1.14). But we are unable to show it.\\

The transversal fluctuation exponent  $\xi(u)$ is related to another important {\em longitudinal exponent} $\chi(u)$.  Let
$$\chi(u)=\sup\{r \geq 0: \mbox{ for  some }C>0, {\bf \sigma}^2(T({\bf 0}, n u))\geq Cn^{2r} \mbox{ for all }n\},$$
where ${\bf \sigma}^2(X)$ is the variance of $X$.
It is conjectured by statistical physicists (see Krug and Spohn (1992)) that
$$\chi(u)= 2\xi(u)-1 \mbox{ for all }u\in {\bf R}^2.\eqno{(1.19)}$$
If $F(0)< p_c$ and infsupp$(F)=0$, or infsupp$(F)=l >0$ but $F(l) < \vec{p}_c$, 
Wehr and Aizenman (1990) (see Theorem 5 in Newman and Piza (1995)) showed that
$$\chi(u) \geq {1-\xi(u)\over 2}.\eqno{(1.20)}$$

When  (1.14) holds, it is known  (see the discussion of case 2 in Newman and Piza (1995)
and Remark 6 in Zhang (2006)) that for all $\theta_p^-\leq \theta \leq \theta_p^+$, we have 
 $${\bf \sigma}^2(({\bf 0}, n (r, \theta)))\leq M.\eqno{(1.21)}$$
This implies that for all $\theta_p^-\leq \theta \leq \theta_p^+$,
$$\chi((r, \theta))=0.\eqno{(1.22)}$$
Both (1.22) and Theorem 2  tell us that conjecture (1.19) cannot be true when $F$ satisfies (1.14).
Furthermore, note that
as we showed in Theorem 3,  $\xi((r, \theta_p^\pm))= 0.5$, so even
Wehr and Aizenman's inequality cannot be true when $F$ satisfies (1.14).

\section{ Proof of Theorem 1.} 
Before the proof of Theorem 1, we would like to introduce a few lemmas. 
We denote 
$$tA=\{(tx, ty): (x,y)\in A\}\mbox{ for }A\subset {\bf R}^2.$$
In particular, when $t >1$ or $t <1$, $tA$ is said to be enlarged or shrunk $t$ times.
\\

{\bf Lemma 1.}  {\em If $F$ satisfies (1.14) with $F(0) < p_c$, 
then for any $v\not\in m{\bf B}_F^\circ$ with $m \geq 1$, there exist $C_i=C_i(F, \delta)$
for $i=1,2$ such that
$${\bf P}\left[ T({\bf 0}, v)\leq m -m^{1/2+\delta}\right]\leq C_1 \exp(-C_2 m^{\delta}).$$
In addition, if $v\in  m{\bf B}_F$, then there exist $C_i=C_i(F, \delta)$
for $i=1,2$ such that}
$${\bf P}\left[ T({\bf 0}, v)\geq m +m^{1/2+\delta}\right]\leq C_1 \exp(-C_2 m^{\delta}).$$

{\bf Proof.} For each $v=(r_v, \theta_v) \in m \partial {\bf B}_F$ with $m\geq 1$, 
we draw a line from the origin  to $v$.
The line has to pass through $\partial {\bf B}_F$ at $v'$. By the definition, we know that $\mu_F(v')\geq 1$
and $m v' =v$. 
By  (1.3) and Theorem 3.1 in Alexander (1997), there exists $C=C(F)$ such that for any $v'\in \partial {\bf B}_F$,
$$ m\leq m \mu_F(v') \leq ET({\bf 0}, v)=ET({\bf 0}, mv')\leq m \mu_F(v') +Cm^{1/2} \log m.$$
If we take $m$ with $m^{1/2+\delta} \geq 2C m^{1/2}\log m$, then
\begin{eqnarray*}
&&{\bf P}\left[T({\bf 0}, v) \leq m-m^{1/2+\delta}\right]\\
&&\leq  {\bf P}\left[T({\bf 0}, mv') \leq ET({\bf 0}, mv')-m^{1/2+\delta}\right]\\
&&\leq  {\bf P}\left[T({\bf 0}, m v') -ET({\bf 0}, mv')\leq  -m^{1/2+\delta}\right].
\end{eqnarray*}
It follows from  Theorem 1, (1.15) in Kesten (1993) that
$$ {\bf P}\left[T({\bf 0}, v) \leq m-m^{1/2+\delta}\right]\leq C_1 \exp(-C_2 m^\delta).\eqno{(2.1)}$$
Lemma 1 follows if $v\in m \partial {\bf B}_F$. For any $v\not\in m{\bf B}_F^\circ$, let 
$\gamma_m$ be an optimal path from the origin to $v$. Then $\gamma_m$ has to meet $m \partial {\bf B}_F$
at $\alpha$. Note that if $\gamma_m$  is not unique, we can select $\gamma_m$ in a unique way. 
Without loss of generality, we assume that $\alpha\in {\bf Z}^2$. Otherwise, we can always select a neighbor
vertex of $\alpha$ in ${\bf Z}^2$ such that it has the same passage time as $T({\bf 0}, \alpha)$.
Therefore, by (2.1)
\begin{eqnarray*}
&&{\bf P}\left[T({\bf 0}, v) \leq m-m^{1/2+\delta}\right]\\
&&\leq {\bf P}\left[T({\bf 0}, \alpha) \leq m-m^{1/2+\delta}\right]\\
&&\leq  \sum_{u\in m\partial {\bf B}_F\cap {\bf Z}^2}{\bf P}\left[T({\bf 0}, v) \leq m-m^{1/2+\delta}, \alpha =u\right]\\
&&\leq C_1m^2\exp(-C_2 m^\delta).
\end{eqnarray*}
The first inequality in Lemma 1  follows. 
Note that if $v\in {\bf B}_F$, then $\mu_F(v) \leq 1$, so
the same proof  can be adapted to show the second probability estimate in Lemma 1. $\Box$\\



Kesten (1986)  proved the following lemma in his Theorem 8.5.\\

{\bf Lemma 2} (Kesten).  {\em If $F(0) < p_c$, for any two vertices $u$ and $v$ with
$\|u-v\|=m$, there exist constants $C_i=C_i(F)$  for $i=1,2,3,4,5,6$ such that}
$${\bf P}\left[T(u,v) \leq C_1 \|u-v\|\right] \leq C_2 \exp(-C_3 m)\mbox{ and }{\bf P}\left[T(u,v) \geq C_4 \|u-v\|\right] \leq C_5 \exp(-C_6 m).$$

With Lemma 2, we also have the following lemma.\\

{\bf Lemma 3.}  {\em For any  vector, let $\gamma_n(u)$ be an optimal path from $(0,0)$ to
$nu$. If $F(0)< p_c$, there exist $C_i=C_i(F, u) $ for $i=1,2,3 $ such that}
$${\bf P}\left[\gamma_n(u)\not \subset [-C_1n,C_1n]^2 \right] \leq C_2 \exp(-C_3 n).$$
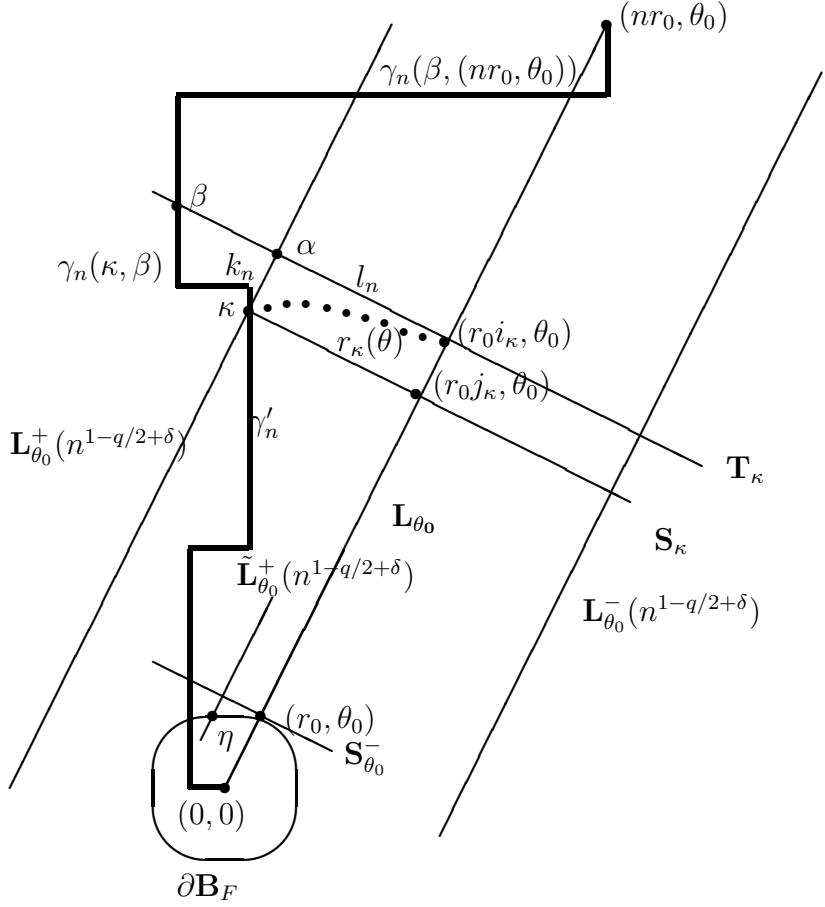
\begin{figure}\label{F:graph}
\begin{center}
\setlength{\unitlength}{0.0125in}%
\begin{picture}(200,230)(67,850)
\thicklines

\put(60,1000){\circle*{4}}
\put(68,1001){\circle*{3}}
\put(76,1003){\circle*{3}}
\put(84,1003){\circle*{3}}
\put(93,1001){\circle*{3}}
\put(101,999){\circle*{3}}
\put(109,997){\circle*{3}}
\put(117,995){\circle*{3}}
\put(125,991){\circle*{3}}
\put(133,989){\circle*{3}}
\put(142,987){\circle*{4}}
\put(130,965){\circle*{4}}

\put(147,987){$(r_0i_\kappa, \theta_0)$}
\put(97,983){$r_\kappa(\theta)$}

\put(137,965){$(r_0j_\kappa, \theta_0)$}
\put(72,1024){\circle*{4}}

\put(80,1024){$\alpha$}
\put(47,998){$\kappa$}
\put(30,1044){\circle*{4}}
\put(35,1044){$\beta$}
\put(115,1097){$\gamma_n(\beta, (nr_0, \theta_0))$}

\put(50,1015){$k_n$}
\put(105,1010){$l_n$}
\put(-20,1015){$\gamma_n(\kappa, \beta)$}
\put(60,950){$\gamma_n'$}
\put(260,930){${\bf T}_\kappa$}
\put(230,900){${\bf S}_\kappa$}

\put(50, 800){\line(1,2){160}}
\put(50,800){\circle*{4}}

\put(-40,940){${\bf L}_{\theta_0}^+(n^{1-q/2+\delta})$}
\put(200,870){${\bf L}_{\theta_0}^-(n^{1-q/2+\delta})$}
\put(-40, 800){\line(1,2){160}}
\put(140, 780){\line(1,2){160}}
\put(50, 800){\line(-1,0){15}}
\put(50, 801){\line(-1,0){15}}
\put(36, 800){\line(0,1){100}}

\put(35, 800){\line(0,1){100}}
\put(35, 900){\line(1,0){25}}
\put(60, 900){\line(0,1){110}}
\put(60, 1010){\line(-1,0){30}}
\put(30, 1010){\line(0,1){80}}
\put(30, 1090){\line(1,0){180}}
\put(210, 1090){\line(0,1){30}}
\put(35, 901){\line(1,0){25}}
\put(61, 900){\line(0,1){110}}
\put(60, 1011){\line(-1,0){30}}
\put(31, 1010){\line(0,1){80}}
\put(30, 1091){\line(1,0){180}}
\put(211, 1090){\line(0,1){30}}

\put(210,1120){\circle*{4}}
\put(215,1120){$(nr_0, \theta_0)$}

\put(60,1000){\line(2,-1){160}}
\put(20,1050){\line(2,-1){230}}

\put(50,800){\oval(60,60)}

\put(65,830){\circle*{4}}

\put(45,830){\circle*{4}}
\put(50, 800){\line(1,2){50}}
\put(40, 820){\line(1,2){30}}
\put(47, 818){$\eta$}
\put(75, 825){$(r_0, \theta_0)$}

\put(55, 885){$\tilde{{\bf L}}_{\theta_0}^+(n^{1-q/2+\delta})$}

\put(20,853){\line(2,-1){75}}
\put(100,810){${{{\bf S}_{\theta_0}^-}}$}
\put(30,785){${(0,0)}$}

\put(30,755){${\partial {\bf B}_F}$}

\put(120,910){${{\bf L_{\theta_0}}}$}

\end{picture}
\end{center}
\vskip 4cm
\caption{{\em The graph shows the optimal path $\gamma_n$ (the boldfaced line in the graph) from the origin to $(nr_0, \theta_0)$.
The path meets ${\bf L}_{\theta_0}^+(n^{1-q/2+\delta})$ at $\kappa$, so ${\cal E}_n^+$ occurs.  
The line ${\bf S}_\kappa$ passes through $\kappa$,  is perpendicular to ${\bf L}_{\theta_0}^+(n^{1-q/2+\delta})$, and meets
${\bf L}_{\theta_0}$ at $(r_0 j_\kappa, \theta_0)$. The line $\tilde{{\bf L}}^+_n$ is the line 
${\bf L}_{\theta_0}^+(n^{1-q/2+\delta})$ as it shrinks toward the center line ${\bf L}_{\theta_0}$  $n$ times. 
$\tilde{{\bf L}}_{\theta_0}^+(n^{1-q/2+\delta})$
meets $\partial {\bf B}_F$ at $\eta$. The curve along $\partial {\bf B}_F$ from $(r_0,\theta_0)$ to
$\eta$ is $s_{\theta_0}(\theta)$ (it does not show in the graph). 
We enlarge $s_{\theta_0}(\theta)$ $n$ times and move it  to $\kappa$ from
$\eta$ to have  $r_\kappa(\theta)$ (the dotted  line in the graph). $r_\kappa(\theta)$ meets ${\bf L}_{\theta_0}$
at $(r_0 i_\kappa, \theta_0)$.  We draw a line ${\bf T}_\kappa$ passing through $(r_0i_\kappa, \theta_0)$.
We divide the optimal path $\gamma_n$
into three pieces: $\gamma_n'$, $\gamma_n(\kappa, \beta)$, and $\gamma_n(\beta, (nr_0, \theta_0))$. 
The distances from $\kappa$ to $\alpha$ and from $\alpha$ to $(r_0i_\kappa, \theta_0)$ are  $k_n$ and $l_n$,
respectively.
By (2.4), $|\gamma_n(\kappa, \beta)|\geq k_n \geq Cn^{1/2+\delta}$.}}
\end{figure}

{\bf Proof of Theorem 1.} Note that if $\max\{\kappa^-(u), \kappa^+(u)\}=0$, then Theorem 1 follows
from Lemma 3 directly.
We suppose that $q=\kappa^-(u)=\max\{\kappa^-(u), \kappa^+(u)\}$
for some $q>0$ and $u=(r_0, \theta_0)\in \partial {\bf B}_F$ with $0\leq \theta_0 \leq \pi/2$
without loss of generality.  Note that
$$r_0=1/\mu_F(r_0, \theta_0).$$
Recall that ${\bf L}_{\theta_0}$ is the line connecting the origin and $(nr_0,\theta_0)$.
In addition, we denote by ${\bf L}_{\theta_0}^+(n^{1-q/2+\delta})$ and ${\bf L}_{\theta_0}^-(n^{1-q/2+\delta})$  the parallel lines  from above and below
$n^{1-q/2+\delta}$ to ${\bf L}_\theta$, respectively (see Fig. 3).
If the
transversal fluctuations
$h_n(r_0, \theta_0) \geq n^{1-l/2+\delta}$, then there exists an optimal path
$\gamma_n$, from the origin to $(nr_0, \theta_0)$, such that it   meets either ${\bf L}_{\theta_0}^+(n^{1-q/2+\delta})$ or
 ${\bf L}_{\theta_0}^-(n^{1-q/2+\delta})$.
Let 
$${\cal E}^+_n = \{\gamma_n \mbox{ meets }{\bf L}^+_n\}\mbox{ and }{\cal E}^-_n = \{\gamma_n \mbox{ meets }{\bf L}^-_n\}.$$ 
Note that it is possible that $\gamma_n$ meets both ${\bf L}_{\theta_0}^+(n^{1-q/2+\delta})$ and ${\bf L}_{\theta_0}^-(n^{1-q/2+\delta})$. We first suppose that ${\cal E}^+_n$ occurs.
To show Theorem 1, we divide the proof
into the following geometric analysis and probability estimate.\\

{\bf Geometric analysis.}
For a curve with  rectangular or  polar coordinates in a certain scale (e.g., inches), if we change the scale (e.g., feet), 
then we may enlarge or shrink the curve. We say the two curves, the original one and the changed one, are 
{\em similar}.

On ${\cal E}^+_n$, we denote  the first intersection point of ${\bf L}_{\theta_0}^+(n^{1-q/2+\delta})$  and $\gamma_n$ by $\kappa$ (see Fig. 3).
Note that $\kappa$ depends on the configurations, so it is a random number.
We draw a line ${\bf S}_\kappa$ passing through $\kappa$, perpendicular to ${\bf L}_{\theta_0}^+(n^{1-q/2+\delta})$, to meet 
${\bf L}_{\theta_0}$ at $(r_0 j_\kappa, \theta_0)$ (see Fig. 3). 
Note that the intersection point  may be on the left side of the origin on ${\bf L}_{\theta_0}$. 
In this case, we just  say  $j_\kappa$ is negative.

If  ${\bf L}_{\theta_0}^+(n^{1-q/2+\delta})$ is shrunk  $n$ times
smaller toward the center line ${\bf L}_{\theta_0}$, we have 
$\tilde{{\bf L}}_{\theta_0}^+(n^{1-q/2+\delta})
=(1/n){\bf L}_{\theta_0}^+(n^{1-q/2+\delta})$ (see Fig. 3).
If $n$ is large, then $\tilde{{\bf L}}_{\theta_0}^+(n^{1-q/2+\delta})$
 will intersect $\partial {\bf B}_F$ at $\eta=(r_\eta, \theta_\eta)$. Let
$s_{\theta_0}(\theta)$ be the piece of the curve of $\partial {\bf B}_F$ from $(r_0, \theta_0)$ to $\eta$ (see Fig. 3).
Now $s_{\theta_0}(\theta)$  is enlarged $n$ times   and we move the enlarged curve $n s_{\theta_0}(\theta)$ parallel to $\kappa$ 
from $\eta$ to have a curve
$r_\kappa(\theta)$. Clearly, $r_\kappa(\theta)$ is similar to $s_{\theta_0}(\theta)$.
We suppose that the intersection point (see Fig. 3)
$$r_\kappa(\theta)\cap {\bf L}_{\theta_0}=(i_\kappa r_0, \theta_0) .\eqno{(2.2)}$$
Note also that the intersection point  may be on the left side of the origin on ${\bf L}_{\theta_0}$. 
In this case, we just  say  $i_\kappa$ is negative.
By the definition of the curve $r_\kappa(\theta)$, we have
$$j_\kappa \leq i_\kappa.$$
 We denote by ${\bf T}_\kappa$ 
the line that passes through $(r_0i_\kappa, \theta_0)$ and is parallel to ${\bf S}_\kappa$ (see Fig. 3).
Note  that ${\bf S}_\kappa$ and ${\bf T}_\kappa $
are parallel, so we  may assume that 
${\bf T}_\kappa$ meets ${{\bf L}}_{\theta_0}^+(n^{1-q/2+\delta})$
 at $\alpha$ and $l_n=\|\alpha-(r_0i_\kappa, \theta_0)\|$ (see Fig. 3).
We also denote (see Fig. 3)
$$k_n=\mbox{{\bf dist}}(\kappa, {\bf T}_k).$$ 
Note that the two curves $s_{\theta_0}(\theta)$ and $r_\kappa(\theta)$ are similar. Note also that
the distance $l_n$ from $\alpha $ to $(r_0i_\kappa, \theta_0)$ (see Fig. 3) satisfies 
$$l_n \geq n^{1-q/2+\delta}.\eqno{(2.3)}$$
So after shrinking $k_n$ and $l_n$  $n$ times smaller, by the assumption in (1.9),
$$ \left({k_n \over n}\right ) \geq C \left({l_n\over n}\right) ^{1/q}\geq C\left ({ n^{1-q/2+\delta}\over n}\right )^{1/q}.$$
Therefore, 
the distance from $\kappa $ to ${\bf T}_\kappa$, $k_n$, is at least 
$$k_n \geq Cn^{[1/2+\delta/q]}\geq Cn^{[1/2+\delta]},\eqno{(2.4)}$$
for some constant $C=C(\theta, r)$.\\

{\bf Probability estimate.}
Now we divide $\gamma_n$ into two pieces, $\gamma_n'(\kappa)$ and 
$\gamma_n''(\kappa)$, where $\gamma_n'(\kappa)$ is the piece from the origin to $\kappa $ along $\gamma_n$, and
$\gamma_n''(\kappa)$ is the rest of the piece from $\kappa$ to $(nr_0, \theta_0)$.
With this definition,
$$T(\gamma_n)=T(\gamma_n'(\kappa))+T(\gamma_n''(\kappa)).\eqno{}$$

Now we divide the following two cases:

(a) $ n^\delta \leq i_\kappa \leq  n-n^\delta $.

(b) $i_\kappa < n^\delta $ or $i_\kappa >n-n^\delta $.

We focus on the easy case, (b), first. Without loss of generality, by  symmetry, 
we may only focus on  $i_\kappa >n-n^\delta $.
By Lemma 3, there exist $C_i=C_i(F)$ for $i=1,2,3,$ such that
$${\bf P}\left[\kappa \not\in C_1[-n,n]^2\right]\leq C_2\exp(-C_3 n).\eqno{(2.5)}$$  
We now estimate the event that
$$T(\gamma_n'(\kappa)) \leq i_\kappa -i_\kappa^{1/2+\delta}.$$
Without loss of generality, we assume that $\kappa$ is an integer vertex. Otherwise, we can always
select an integer vertex next to $\kappa$ with a distance to $\kappa$ less than 1. Hence, by the definition of $\kappa$ and (2.5),
\begin{eqnarray*}
&&{\bf P}\left[T(\gamma_n'(\kappa)) \leq i_\kappa -i_\kappa^{1/2+\delta},{\cal E}_n^+, \mbox{ (b) occurs}\right]\\
&\leq &\sum_{v\in [-Cn,Cn]^2, v\not\in i_v {\bf B}_F^\circ}{\bf P}\left[T(\gamma_n'(v)) \leq i_v -i_v^{1/2+\delta}, \kappa=v\right ]+C_1\exp(-C_2 n).
\end{eqnarray*}
For a fixed $v$, as we defined, $v$ is out of $i_v {\bf B}_F^\circ $, so by Lemma 1 and above inequality,
$${\bf P}\left[T(\gamma_n'(\kappa)) \leq i_\kappa  -i_\kappa^{1/2+\delta},{\cal E}_n^+, \mbox{ (b) occurs}\right]\leq C_1 n^2\exp(-C_2 i_v^{\delta})\leq C_3\exp(-C_4 n^\delta).\eqno{(2.6)}$$
Note that, for small $\delta$, by the assumption of (b), 
$$i_\kappa -(i_\kappa)^{1/2+\delta/2}\geq (n- n^\delta)-(n- n^\delta)^{1/2+\delta/2},$$ 
so by (2.6) we have
$${\bf P}\left[T(\gamma_n'(\kappa)) \leq n -n^\delta-(n-n^\delta)^{1/2+\delta/2},{\cal E}_n^+, \mbox{ (b) occurs}\right]\leq C_1\exp(-C_2 n^{\delta}).\eqno{(2.7)}$$
Now we estimate $T(\gamma''(\kappa))$. Note that $\kappa$ is on ${{\bf L}}_{\theta_0}^+(n^{1-q/2+\delta})$,
 so the distance from $\kappa$ to
 $(nr_0, \theta_0)$ is at least $n^{1-q/2+\delta}$. By Lemma 2, there exist $C_i(\delta, F)$ for $i=1,2,3,4,5,6,7$
such that
\begin{eqnarray*}
&&{\bf P}\left [T(\gamma_n''(\kappa)) \leq C_1n^{1-q/2+\delta}, {\cal E}_n^+,  \mbox{(b) occurs}\right]\\
&&\leq \sum_{v\in [-Cn,Cn]^2}{\bf P}\left[T(\gamma_n''(v)) \leq C_1n^{1-q/2+\delta}, \kappa =v,{\cal E}_n^+,  \mbox{(b) occurs}\right] +C_2\exp(-C_3 n)\\
&& \leq C_4 n^2 \exp(-C_5 n^{\delta})+C_2\exp(-C_3 n)\\
&&\leq C_6 \exp(-C_7 n^{\delta}).\hskip 4in (2.8)
\end{eqnarray*}

For small $\delta$
and $C_1$ defined in (2.8), we take large $n$ such that
$$n^\delta+(n- n^\delta)^{1/2+\delta/2}\leq C_1n^{1/2+\delta}/2.$$
Note that $n^{1-q/2+\delta} \geq n^{1/2+\delta}$, so with (2.7), (2.8), and Lemma 1,
\begin{eqnarray*}
&&{\bf P}\left[ h_n(u) \geq n^{1-(q/2)+\delta} , \mbox{ (b) occurs}, {\cal E}_n^+\right]\\
&&\leq  {\bf P}\left[ h_n(u) \geq n^{1-(q/2)+\delta} , {\cal E}_n^+, T(\gamma_n'(\kappa)) \geq n  -C_1n^{1/2+\delta}/2,
T(\gamma_n''(\kappa) )\geq C_1n^{1-(q/2)+\delta}\right]\\
&&\hskip 1cm +C_2\exp(-C_3 n^{\delta})\\
&&\leq {\bf P}\left[T(\gamma_n) \geq n +C_1n^{1/2+\delta}/2\right]+C_2\exp(-C_3 n^{\delta})\\
&&\leq C_4\exp(-C_5 n^{\delta}).\hskip 4.2in (2.9)
\end{eqnarray*}

Now we focus on the event that 
${\cal E}_n^+$ and (a) occurs. 
Under the assumption of (a), ${\bf T}_\kappa$ is between the origin and $(nr_0, \theta_0)$, so
$\gamma_n$ meets ${\bf T}_\kappa$ at $\beta$ (see Fig. 3).
Note that path $\gamma_n$  reaches to $\kappa$, to $\beta\in {\bf T}_\kappa$, and then to
$(nr_0, \theta_0)$. Recall that path $\gamma_n'(\kappa)$ is from the origin to $\kappa$.
In addition, we divide $\gamma''(n)$ into two pieces: 
$\gamma_n(\kappa, \beta)$ from $\kappa$ to $\beta\in {\bf T}_\kappa$ along $\gamma_n$,
and $\gamma_n(\beta, (nr_0, \theta_0))$ from $\beta $ to $(nr_0, \theta_0)$. We have
$$T(\gamma''(n))= T(\gamma_n(\kappa, \beta))+ T(\gamma_n(\beta, (n r_0,\theta_0))).$$
We want to remark that it is possible that $\gamma_n$ reaches ${\bf T}_\kappa$
first, then moves on to $\kappa$. By (2.4), 
$$|\gamma_n(\kappa, \beta)|\geq n^{1/2+\delta}.\eqno{(2.10)}$$
By Lemma 3 and (2.10), to fix both $\kappa$ and $\beta$, we have
\begin{eqnarray*}
&&{\bf P}\left[T(\gamma_n(\kappa, \beta))\leq n^{1/2+2\delta/3}, {\cal E}^+_n, \mbox{ (a) occurs}\right]\\
&\leq &\sum_{v,u\in [-Cn, Cn]^2, \|u-v\| \geq n^{1/2+\delta}} {\bf P}\left[T(u, v)\leq n^{1/2+2\delta/3}\right]+C_1\exp(-C_2 n).
\hskip 1in (2.11)
\end{eqnarray*}
By Lemma 2 and (2.11),
$${\bf P}\left[T(\gamma_n(\kappa, \beta))\leq n^{1/2+2\delta/3},{\cal E}^+_n, \mbox{ (a) occurs}\right]\leq C_1n^4\exp(-C_2 n^{1/2})\leq C_3\exp(-C_4 n^{1/2}).\eqno{(2.12)}$$

Now we work on $\gamma_n'(\kappa)$. Note that $i_\kappa \geq n^\delta$, so by the same estimate of (2.6), we have 
$${\bf P}\left[T(\gamma_n'(\kappa)) \leq i_\kappa -i_\kappa^{1/2+\delta/2}\right]\leq C_1\exp(-C_2 n^{\delta/2}).\eqno{(2.13)}$$

Finally, we focus on  $\gamma_n(\beta, (nr_0,\theta_0))$. We reconsider
$(nr_0, \theta_0)$ as the new origin, and   ${\bf B}_F(nr_0, \theta_0)$ as the shape by moving 
${\bf B}_F$ from the origin to $(nr_0, \theta_0)$. We enlarge ${\bf B}_F(nr_0, \theta_0)$ 
$n-i_\kappa $ times to have $(n- i_\kappa ){\bf B}_F(nr_0, \theta_0)$. Note that
 $(i_\kappa r_0, \theta_0)$ is the intersection of $(n- i_\kappa )\partial {\bf B}_F(nr_0, \theta_0)$
and ${\bf L}_{\theta_0}$.
By  symmetry and the convexity of $(n- i_\kappa ){\bf B}_F(nr_0, \theta_0)$, $(n- i_\kappa ){\bf B}_F(nr_0, \theta_0)$ stays
above the line ${\bf T}_\kappa$  (see Fig. 3). Therefore, $\beta$
 is out side the interior of $(n- i_\kappa ){\bf B}_F(nr_0, \theta_0)$ (see Fig. 3).
By this observation, note that $i_\kappa \leq n-n^\delta$, so by Lemma 3, translation invariance, and Lemma 1, 
\begin{eqnarray*}
&&{\bf P}\left[T(\gamma_n(\beta, (nr_0, \theta_0)) \leq ( n-i_\kappa ) -(n -i_\kappa )^{1/2+\delta/2}\right]\\
&&\leq \sum_{\scriptstyle{v, u\in [-Cn, Cn]^2,}\atop{u\not\in (n- i_v ){{\bf B}}_F^\circ(nr_0, \theta_0), n-i_v \geq n^\delta}} {\bf P}\left[T(\gamma_n(u, (nr_0, \theta_0)) \leq (n-i_v ) -(n-i_v )^{1/2+\delta/2}, \kappa=v, 
\beta=u\right]\\
&& \hskip 2cm +C_1\exp(-C_2 n)\\
&& =\sum_{u\not\in (n- i_v ){{\bf B}}_F^\circ, n-i_v \geq n^\delta}
 {\bf P}\left[T({\bf 0}, u) \leq (n-i_v) -(n-i_v)^{1/2+\delta/2}\right]
+C_1\exp(-C_2 n)\\
&&\leq C_2 n^4 \exp(-C_3 n^{\delta/2}).\hskip 4.3in (2.14)
\end{eqnarray*}
Note that 
$$T(\gamma_n)=T(\gamma_n'(\kappa) )+T(\gamma_n(\kappa, \beta))+T(\gamma_n(\beta, (nr_0, \theta_0))).$$
Note also that $0< i_\kappa < n$, so if
$$T(\gamma_n'(\kappa)) \geq i_\kappa -i_\kappa^{1/2+\delta/2}, T(\gamma_n(\kappa, \beta))\geq n^{1/2+2\delta/3},
\mbox{ and } T(\beta, (nr_0, \theta_0)) \geq (n -i_\kappa)-(n-i_\kappa)^{1/2+\delta/2},$$
then there exists $C$ such that for all large $n$,
$$T(\gamma_n)=T(\gamma_n'(\kappa) )+T(\gamma_n(\kappa, \beta))+T(\gamma_n(\beta, (nr, \theta)))
\geq n+ C n^{1/2+2\delta/3}.$$
Therefore, by (2.12), (2.13), (2.14), and Lemma 1, we have
\begin{eqnarray*}
&&{\bf P}\left[ h_n(u) \geq n^{1-(q/2)+\delta} , \mbox{ (a) occurs}, {\cal E}_n^+\right]\\
&&\leq  {\bf P}[ h_n(u) \geq n^{1-(q/2)+\delta} , {\cal E}_n^+, T(\gamma_n'(\kappa)) \geq i_\kappa -i^{1/2+\delta/2}, T(\gamma_n(\kappa, \beta))\geq n^{1/2+2\delta/3},\\
&&\hskip 2cm T(\gamma_n(\beta, (nr, \theta))) \geq (n -i_\kappa)-(n-i_v)^{1/2+\delta}]+C_1\exp(-C_2 n^{\delta/2})\\
&&\leq {\bf P}\left[T(\gamma_n) \geq n+Cn^{1/2+2\delta/3}\right]+C_1\exp(-C_2 n^{\delta/2})\\
&&\leq C_3\exp(-C_4 n^{{\delta}/2}).\hskip 4.3in (2.15)
\end{eqnarray*}
With (2.9) and (2.15) together, we have
$${\bf P}\left[ h_n(u) \geq n^{1-(q/2)+\delta} , {\cal E}_n^+\right]\leq C_1 \exp(-C_2 n^{\delta/2}).\eqno{(2.16)}$$

Now we focus on the case that ${\cal E}^-_n$ occurs. Note that if we rotate the graph $180^\circ$ around the point
$(nr_0, \theta_0)$, then by  symmetry, translation invariance, and the same estimate of (2.16), we can show
$${\bf P}\left[ h_n(u) \geq n^{1-(1/2l)+\delta} , {\cal E}_n^-\right]\leq C_1 \exp(-C_2 n^{{\delta}/2}).\eqno{(2.17)}$$
Since
$${\bf P}\left[ h_n(u) \geq n^{1-(1/2l)+\delta} \right] \leq {\bf P}\left[ h_n(u) \geq n^{1-(1/2l)+\delta},{\cal E}_n^+\right]
+{\bf P}\left[ h_n(u) \geq n^{1-(1/2l)+\delta},{\cal E}_n^-\right],$$
therefore, Theorem 1 follows from (2.16) and (2.17).\\

\section{ Proofs of Theorems 2 and 3.}
In this section, we will show Theorems 2 and 3. Before showing them, we will introduce a few notations and a lemma.
We denote  the southwest diagonal line passing through $u$  by 
$${\bf D}_u=\{u+(x,y):(x,y):x+y=0\}.$$ 
We also denote by
$${\bf D}_u^+=\{u+(x,y): x\leq 0, x+y=0\}\mbox{ and }{\bf D}_u^-=\{u+(x,y): x\geq 0, x+y=0\}$$
 the upper and lower lines (from $u$ to $\infty$ and
to $-\infty$). Clearly,
$${\bf D}_u={\bf D}_u^+\cup {\bf D}_u^-.$$
Recall from section 2 that  ${\bf L}_{\theta^-_p}$ is defined as the line passing through the origin and  $(r, \theta_p^-)$,
and 
${\bf L}^\pm_{\theta^-_p}(m)$ is the lines parallel  to ${\bf L}_{\theta^-_p}$ from above and below $m$ units. Let
${\bf S}_{\theta^-_p}(m)$ be the strip from ${\bf L}_{\theta_p^-}^-(m)$ to ${\bf L}_{\theta_p^-}^+(m)$.
With this definition, we show the following lemma.\\

{\bf Lemma 4.}{\em If $p=F(0) > \vec{p}_c$,  for any $\delta >0$, 
there exist $C_i=C_i(F, \delta)$ for $i=1,2$ such that}
$${\bf P}\left[{\bf D}_0 \stackrel{\gamma_n}{\rightarrow}{\bf D}_{(n,n)} \mbox{ with }\gamma_n\subset {\bf S}_{\theta^-_p}(n^{1/2+\delta})\right]\geq
1-C_1\exp(-C_2 n^\delta ).\eqno{(3.1)}$$

Before proving Lemma 4, we require some terminology
(see Kuczek (1989)). We use the graph ${\cal L}$ in section 1.
 Given vertices $u$ and $v$ in ${\cal L}$, we say $u\rightarrow v$ if
there is a sequence $v_0=u, v_1,\cdots, v_m=v$ of points of ${\cal
L}$ with $v_i:=(x_i,y_i)$ and $v_{i+1}:=(x_{i+1}, y_i + 1)$ for
$0\leq i\leq m-1$ such that $v_i$ and $v_{i+1}$ are connected by
an edge with weight 1. Thus, $u\rightarrow v$ if there is a
sequence of oriented edges, each with weight 1, joining $u$ to
$v$. For two sets $A$ and $B$, we denote $A\rightarrow B$ if
there are $u\in A$ and $v\in B$ such that $u\rightarrow v$. 
 For
$A\subset {\bf Z}$, let
$$\xi_n^A :=\{x: (x,n) \in {\cal L}
\mbox{ and}  \,\,\, \exists \,\,\, x'\in A \mbox{ such that }
(x',0) \rightarrow (x,n) \mbox{ for } n>0\}.$$

As in Kuczek (1989), denote the event that there exists an
infinite oriented path of weight-1 edges starting from $(x, y)$ by
$\Omega_\infty^{(x,y)}$. We let $\xi_0':=\xi_0^{(0,0)}:= (0,0)$
and
set \\
$$\xi_1':=\left \{ \begin{array}{ll}
                     \xi_1^{(0,0)} &\mbox{ if $\xi_1^{(0,0)}\neq \emptyset,$}\\
                     \{1\} &\mbox{ otherwise,}
                   \end{array}
                   \right.
                    $$
\\
and define inductively, for all $n = 1,2,...$
\\
$$\xi_{n+1}':=\left \{ \begin{array}{ll}
                     \{ x: (x,n+1) \in {\cal L} \mbox{ and }
                      (y,n) \rightarrow (x, n+1) \mbox{ for some }y\in \xi_n'\}
                 &\mbox{ if this set is non-empty}\\
                     \{n+1\} &\mbox{ otherwise.}
                   \end{array}
                   \right.
                    $$
\\ We have suppressed the dependence of $\xi_n'$ on $p=F(1)$ for
notational convenience. Note that $\xi_n'$  is a subset of the
integers between $-n$ and $n$.  Let
$$
r_n'(p):=\sup \xi'_n.
$$
On $\{\xi_n^{(0,0)}\neq \emptyset\}$, we have equivalence between
$r_n'(p)$ and the right-hand edge $r_n(p).$   A vertex $(x, n)\in
{\cal L}$ is said to be a {\em percolation point} if and only if
the event  $\Omega^{(x,n)}_\infty$ occurs.  Let
\begin{eqnarray*}
&&T_1:=\inf\{n\geq 1: (r_n',n) \mbox{ is a percolation point}\},
\,\,\,\, \mbox{ }\\
&& T_2:=\inf\{n\geq T_1+1: (r_n',n) \mbox{ is a percolation point}\}, \,\,\,\, \\
&&\cdots\\
&& T_m:=\inf\{n\geq T_{m-1}+1: (r_n',n) \mbox{ is a percolation
point}\},
\end{eqnarray*}
where we make the convention that $\inf \emptyset =\infty$. Define
$$\tau_1 :=T_1, \ \tau_2 := T_2-T_1,\cdots, \ \tau_m := T_m-T_{m-1},\eqno{(3.2)}$$
where $\tau_i:=0$ if $T_i$ and $T_{i-1}$ are infinite. (Note that
$T_i$ and $T_{i-1}$ are finite with probability 1.)
 Also define
$$ X_1:= r'_{T_1}, \ \ X_2:=r'_{T_2}- r'_{T_1},\cdots, \ \ X_{m}:=r'_{T_{m}}-r'_{T_{m-1}},\eqno{(3.3)}$$
where $X_i := 0$ if $T_i=\infty$ and $T_{i-1}=\infty$. The points
$\{(r'_{T_i}, T_i ) \}$ are called {\em break points} 
since they break the behavior of the right-hand edge into i.i.d.
pieces when the origin is a percolation point. Kuczek (Theorem on
page 1324 (1989)) proved that conditional on
$\Omega^{(0,0)}_\infty$, $\{(X_i, \tau_i)\}$ are i.i.d. with all
moments.  Moreover, for all $q \in (\vec{p}_c, 1]$
 there exists a constant $C_2:=C_2(q)$  such that for all $p \in
 [q,1]$ and all $t \geq 1$,
$$
{\bf P}_p [ X_1 \geq t ] \leq {\bf P}_p [ \tau_1 \geq t ] \leq {\bf P}_p [\xi_{t-1}^{(1,1)} \neq
\emptyset, (1,1) \not\rightarrow \infty] \leq 
C_1\exp(-C_2t), \eqno{(3.4)}
$$
where the last inequality is as in Durrett's sect. 12 (1984).


Now we  use of the following
probability measure on $\Omega_\infty^{(0,0)}$:
$$\bar {\bf P}[\  \cdot \ ]:=
{\bf P} [\ \cdot \, \ \, | \ \ \Omega^{(0,0)}_\infty ]. \eqno{(3.5)}
$$
Let $ \bar{\bf E}$ denote the expected value with respect to
$\bar {\bf P}$.  With these definitions, we begin to show Lemma 4.\\



{\bf Proof of Lemma 4.} We rotate the ${\bf Z}^2$ lattice $45^\circ$ with correct dilation to have ${\cal L}$.
After the rotation, ${\bf D}_{\bf 0}$ and ${\bf D}_{(n,n)}$ will be the $X$-axis and the line $\{y=2n\}$, respectively. 
In addition, after the rotation, 
${\bf L}_{\theta_p^-}^\pm (n^{1/2+\delta})$ and ${\bf S}_{\theta_p^-}(m)$ 
will be denoted by $\bar{\bf L}_{\theta_p^-}^\pm (n^{1/2+\delta})$ and $\bar{\bf S}_{\theta_p^-}(m)$, 
respectively. Thus, to prove Lemma 4, we only need to show
$${\bf P}\left[\exists \mbox{ a $1$-path $\gamma_n$, $\{y=0\}\stackrel{\gamma_n}{\rightarrow} \{y=2n\}$ with $\gamma_n\subset \bar{{\bf S}}(n^{1/2+\delta})$}\right]
\geq 1-C_1\exp(-C_2 n^\delta ).\eqno{(3.6)}$$
We will first  show that
$$\bar{{\bf P}}\left[ r_{2n} \cap \bar{{\bf L}}_{\theta_p^-}^- (n^{1/2+\delta}) \neq \emptyset\right]\leq C_1 \exp(-C_2n^\delta)\eqno{(3.7)}$$
and $$\bar{{\bf P}}\left[ r_{2n} \cap \bar{{\bf L}}_{\theta_p^-}^+ (n^{1/2+\delta}) \neq \emptyset\right]\leq C_1 \exp(-C_2n^\delta).\eqno{(3.8)}$$
We only show (3.7). The same proof can be adapted to show (3.8).\\

Note that, on $\Omega_\infty^{(0,0)}$,  $r_{2n}$ is  the $1$-path from the origin to $(r_{T_1}=X_1, T_1)$, from $(r_{T_1}, T_1)$ to $(r_{T_2}, T_2)$,
$\cdots$, from $(r_{T_{k-1}},T_{k-1})$ to $(r_{T_k}, T_k)$.
On $\{r_{2n} \cap \bar{{\bf L}}_{\theta_p^-}^- (n^{1/2+\delta}) \neq \emptyset\}$, let
$$s_n=\inf \{m\geq 1: \mbox{the $1$-path along $r_{2n}$ from $(r_{T_m}, T_m)$ to $(r_{T_{m+1}}, T_{m+1})$}\}\cap \bar{{\bf L}}_{\theta_p^-}^+ (n^{1/2+\delta})\neq \emptyset\}.$$
On $r_{2n}\cap \bar{{\bf L}}_{\theta_p^-}^+ (n^{1/2+\delta})\neq \emptyset$,
if $s_n=i$, let the first intersection of 
$r_{2n}$ and $\bar{{\bf L}}_{\theta_p^-}^+ (n^{1/2+\delta})$ be $(x_I, y_I).$
By the definition, note that $|X_i|\leq \tau_i$, so
$$ X_1+\cdots +X_{i-1} =r_{T_{i-1}}\leq x_I\leq   X_1+\cdots +X_i+\tau_i.\eqno{(3.9)}$$
Note that $\bar{{\bf L}}_{\theta_p^-}^+(n^{1/2+\delta})$ is the line with the equation
$$y=\alpha_p^{-1} x+\sqrt{2}n^{1/2+\delta},$$
so by (3.9),
$$\alpha_p (\tau_1+\cdots +\tau_{i-1})+\sqrt{2} n^{1/2+\delta} \leq x_I \leq X_1+X_2+\cdots X_{i-1}+\tau_i.\eqno{(3.10)}$$
With these observations,
\begin{eqnarray*}
&&\bar{{\bf P}}\left [r_{2n} \cap \bar{{\bf L}}_{\theta_p^-}^- (n^{1/2+\delta})\neq \emptyset\right]\\
&&\leq \sum_{i=1}\bar{{\bf P}}\left[\alpha_p (\tau_1+\cdots +\tau_{i-1})+n^{1/2+\delta} \leq x_I \leq X_1+\cdots X_{i-1}+\tau_i,s_n=i\right]\\
&&= \sum_{i\leq n^{\delta}}\bar{\bf {P}}\left[\alpha_p (\tau_1+\cdots +\tau_{i-1})+n^{1/2+\delta} \leq x_I \leq X_1+\cdots X_{i-1}+\tau_i,s_n=i\right]\\
&&+\sum_{n\geq i> n^{\delta}}\bar{{\bf P}}\left[\alpha_p (\tau_1+\cdots +\tau_{i-1})+n^{1/2+\delta} \leq x_I \leq X_1+\cdots X_{i-1}+\tau_i,s_n=i\right].
\end{eqnarray*}
Denote the above two sums in the last equations by $I$ and $II$.  Note that $|X_j|\leq \tau_j$, so
\begin{eqnarray*}
&&I:=\sum_{i\leq n^{\delta}}\bar{{\bf P}}\left[\alpha_p (\tau_1+\cdots +\tau_{i-1})+n^{1/2+\delta} \leq x_I \leq X_1+\cdots X_{i-1}+\tau_i,s_n=i\right]\\
&&\leq \sum_{i\leq n^{\delta}} \bar{{\bf P}}\left[n^{1/2+\delta} \leq \tau_1+\cdots +\tau_i\right]
\end{eqnarray*}
Since $\{\tau_i\}$ are i.i.d. with an exponential smaller tail, by a Markov's inequality, we have
$$I\leq C_1\exp(-C_2 n^{1/2}).\eqno{(3.11)}$$
Let estimate $II$. 
\begin{eqnarray*}
&&II:=\sum_{n\geq i> n^{\delta}}\bar{{\bf P}}\left[\alpha_p (\tau_1+\cdots +\tau_{i-1})+n^{1/2+\delta} \leq x_I\leq X_1+\cdots X_{i-1}+\tau_i,s_n=i\right]\\
&&\leq \sum_{n\geq i> n^{\delta}}\bar{{\bf P}}\left[\alpha_p (\tau_1+\cdots +\tau_{i-1})+n^{1/2+\delta} \leq X_1+\cdots X_{i-1}+\tau_i,s_n=i\right]\\
&&\leq \sum_{n\geq i> n^{\delta} }\bar{{\bf P}}\left[ n^{1/2+\delta}/2 \leq  (X_1-\alpha_p \tau_1)+\cdots + (X_{i-1}-\alpha_p \tau_{i-1}) , s_n=i\right]\\
&&+\sum_{n\geq i> n^{\delta} }\bar{{\bf P}}\left[n^{1/2+\delta}/2 \leq  \tau_i, s_n=i\right].\hskip 3.2in (3.12)
\end{eqnarray*}
By (3.4),  the last sum in (3.12)
$$\sum_{n\geq i> n^{\delta} }\bar{{\bf P}}\left[n^{1/2+\delta}/2 \leq  \tau_i, s_n=i\right]\leq 
n \bar{{\bf P}}\left[ n^{1/2+\delta}/2 \leq X_1\right]\leq C_1\exp(-C_2 n^{1/2}).\eqno{(3.13)}$$
Let us estimate the first sum on the right side of (3.12). By Lemma 1 of Zhang (2004), we know that $\bar{E}(X_i-\alpha_p \tau_i)=0$.
In addition, as we mentioned, $X_i-\alpha_p \tau_i$ and $X_j-\alpha_p \tau_j$ are independent with
 exponentially smaller tails, so by a standard large deviation result,  we have
$$\sum_{n\geq i\geq n^{\delta} }\bar{{\bf P}}\left[ n^{1/2+\delta}/2 \leq  (X_1-\alpha_p \tau_1)+\cdots + (X_{i-1}-\alpha_p \tau_{i-1}) , s_n=i\right]\leq C_1 \exp(-C_2 n^\delta).\eqno{(3.14)}$$
Therefore, (3.7) follows from (3.11), (3.13) and (3.14).

Now we use (3.7) and (3.8) to show (3.6). 
By using the inequality in section 12 of Durrett (1984), we have
$${\bf P}\left[\xi^{[-n^{\delta}, 0]}_n\neq \emptyset \right]\geq 1-C_1\exp(-C_2 n^{\delta}).\eqno{(3.15)}
$$
On $\xi^{[-n^{\delta}, 0]}_n\neq \emptyset $, let 
$$\tau=\max\left\{t\leq 0: \Omega_\infty^{(t,0)} \mbox{ occurs}\right\} \mbox{ and } r_n(k)=\sup \xi^{(-\infty,k]}_n .$$
With this definition, by (3.15), translation invariance, (3.7), and (3.8),
\begin{eqnarray*}
&&1-C_1\exp(-C_2 n^{\delta})\\
&&\leq {\bf P}\left[\xi^{[-n^{\delta}, 0]}_n\neq \emptyset \right]\\
&&= \sum_{-n^{\delta}\leq k\leq 0} {\bf P}\left[ \Omega_\infty^{(k,0)},\tau=k\right]\\
&&\leq   \sum_{-n^{\delta}\leq k\leq 0} {\bf P}\left[ \Omega_\infty^{(k,0)}, r_{2n}(k)\cap \{(k,0)+\bar{{\bf S}}_{\theta_p^-}(n^{1/2+\delta}/2)\}= \emptyset ,\tau=k\right]\\
&&\hskip 1cm + n^\delta {\bf P}\left [ \Omega_\infty^{(0,0)}\right] \left(\bar{\bf P}\left[r_{2n} \cap  \bar{{\bf L}}_{\theta_p^-}^+ (n^{1/2+\delta}/2)\neq \emptyset\right]+\bar{\bf P}\left[r_{2n} \cap  \bar{{\bf L}}_{\theta_p^-}^-(n^{1/2+\delta}/2)\neq \emptyset\right]\right)\\
&&\leq  \sum_{-n^{\delta}\leq k\leq 0} {\bf P}\left[ \Omega_\infty^{(k,0)}, r_{2n}(k)\cap \{(k,0)+\bar{{\bf L}}_{\theta_p^-}^\pm (n^{1/2+\delta})\}= \emptyset ,\tau=k \right]+C_1\exp(-C_2 n^\delta).\hskip 0.5cm (3.16)
\end{eqnarray*}
Clearly, for large $n$,
\begin{eqnarray*}
&&\bigcup_{-n^\delta \leq k\leq 0}\left\{\Omega_\infty^{(k,0)}, r_{2n}(k)\cap \left\{(k,0)+\bar{{\bf L}}_{\theta_p^-}^\pm (n^{1/2+\delta}/2)\right\}= \emptyset, \tau=k \right\}\\
&\subset &\left\{\mbox{$\exists$ a $1$-path $\gamma_n$, $\{y=0\}\stackrel{\gamma_n}{\rightarrow} \{y=2n\}$ with $\gamma_n\subset \bar{{\bf S}}(n^{1/2+\delta})$}\right\}.\hskip 3.2cm (3.17)
\end{eqnarray*}
Therefore, (3.6) follows, so does Lemma 4. 
 $\Box$\\

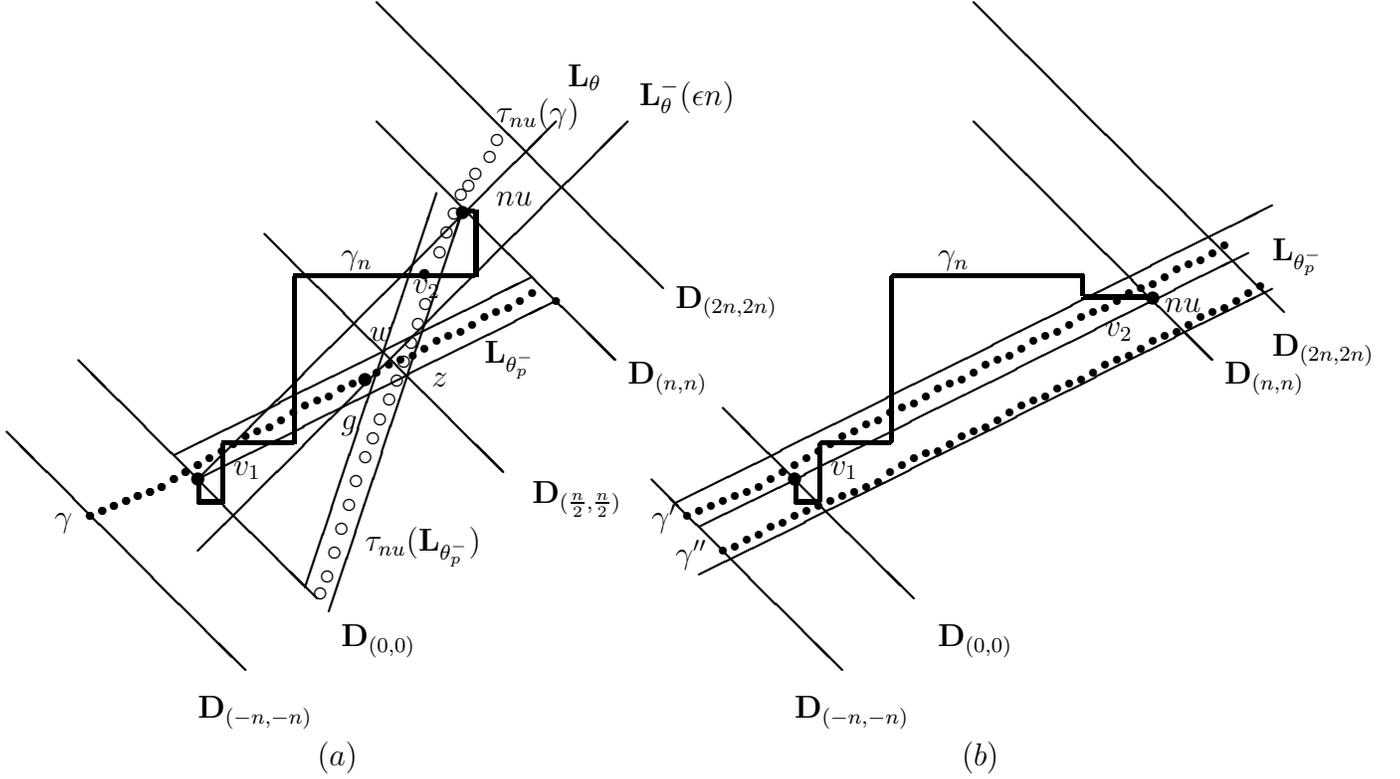
\begin{figure}\label{F:graph}
\begin{center}
\setlength{\unitlength}{0.0125in}%
\begin{picture}(200,230)(107,700)
\thicklines
\put(150,875){\circle*{3}}
\put(0, 800){\line(1,1){150}}
\put(0, 770){\line(1,1){180}}

\put(0, 800){\line(2,1){150}}
\put(60, 730){${\bf D}_{(0,0)}$}
\put(70, 770){${\tau}_{nu}({\bf L}_{\theta_p^-})$}

\put(-10, 810){\line(2,1){150}}
\put(0, 700){${\bf D}_{(-n,-n)}$}

\put(0, 800){\line(1,-1){50}}
\put(0, 800){\line(-1,1){50}}
\put(-30, 770){\line(1,-1){50}}
\put(-30, 770){\line(-1,1){50}}
\put(75,950){\line(1,-1){100}}
\put(75,1000){\line(1,-1){120}}
\put(28,903){\line(1,-1){100}}
\put(140, 790){${\bf D}_{({n\over2},{n\over2})}$}


\put(45, 755){\line(1,3){55}}
\put(55, 745){\line(1,3){56}}

\put(90,878){$v_2$}
\put(15,802){$v_1$}
\put(265,802){$v_1$}
\put(72,858){$w$}

\put(125,915){$nu$}
\put(180, 840){${\bf D}_{(n,n)}$}
\put(200, 872){${\bf D}_{(2n,2n)}$}
\put(400, 876){\circle*{5}}
\put(405, 870){$nu$}
\put(380, 860){$v_2$}

\put(430, 840){${\bf D}_{(n,n)}$}
\put(450, 852){${\bf D}_{(2n,2n)}$}

\put(111, 912){\circle*{6}}
\put(155,965){${\bf L}_\theta$}
\put(125,950){$\tau_{nu}(\gamma)$}

\put(185,957){${\bf L}_\theta^-(\epsilon n)$}

\put(122,939){$\circ$}
\put(119,932){$\circ$}
\put(113,925){$\circ$}
\put(110,920){$\circ$}

\put(107,916){$\circ$}
\put(104,908){$\circ$}
\put(101,900){$\circ$}
\put(98,892){$\circ$}
\put(95,886){\circle*{4}}
\put(92,870){$\circ$}
\put(89,862){$\circ$}
\put(86,854){$\circ$}
\put(83,846){$\circ$}
\put(80,838){$\circ$}
\put(77,830){$\circ$}
\put(73,822){$\circ$}
\put(70,814){$\circ$}
\put(67,806){$\circ$}
\put(63,796){$\circ$}
\put(60,786){$\circ$}
\put(57,776){$\circ$}
\put(54,766){$\circ$}
\put(51,756){$\circ$}
\put(48,749){$\circ$}

\put(-60,780){$\gamma$}
\put(190,780){$\gamma'$}
\put(200,765){$\gamma''$}

\put(-45,785){\circle*{3}}
\put(-40,787){\circle*{3}}
\put(-35,789){\circle*{3}}
\put(-30,791){\circle*{3}}
\put(-25,793){\circle*{3}}
\put(-20,795){\circle*{3}}
\put(-15,797){\circle*{3}}
\put(-10,800){\circle*{3}}
\put(-5,803){\circle*{3}}
\put(0,806){\circle*{3}}
\put(-45,785){\circle*{3}}
\put(-40,787){\circle*{3}}
\put(-35,789){\circle*{3}}
\put(-30,791){\circle*{3}}
\put(-25,793){\circle*{3}}
\put(-20,795){\circle*{3}}
\put(-15,797){\circle*{3}}
\put(-10,800){\circle*{3}}
\put(-5,803){\circle*{3}}
\put(0,806){\circle*{3}}
\put(5,809){\circle*{3}}
\put(10,812){\circle*{3}}
\put(15,815){\circle*{3}}
\put(20,818){\circle*{3}}
\put(25,820){\circle*{3}}
\put(30,822){\circle*{3}}
\put(35,825){\circle*{3}}
\put(40,828){\circle*{3}}
\put(45,831){\circle*{3}}
\put(50,833){\circle*{3}}
\put(55,835){\circle*{3}}
\put(60,838){\circle*{3}}
\put(65,840){\circle*{3}}
\put(70,842){\circle*{5}}
\put(75,845){\circle*{3}}
\put(80,848){\circle*{3}}
\put(85,850){\circle*{3}}
\put(90,852){\circle*{3}}
\put(95,855){\circle*{3}}
\put(100,858){\circle*{3}}
\put(105,860){\circle*{3}}
\put(110,862){\circle*{3}}
\put(115,865){\circle*{3}}
\put(120,868){\circle*{3}}
\put(125,870){\circle*{3}}
\put(130,872){\circle*{3}}
\put(135,875){\circle*{3}}
\put(140,878){\circle*{3}}
\put(60,821){$g$}

\put(0, 800){\line(0,-1){10}}
\put(1, 800){\line(0,-1){10}}

\put(0, 790){\line(1,0){10}}
\put(0, 791){\line(1,0){10}}

\put(10, 790){\line(0,1){25}}
\put(11, 790){\line(0,1){25}}

\put(10, 815){\line(1,0){30}}
\put(10, 816){\line(1,0){30}}

\put(40, 815){\line(0,1){30}}
\put(41, 815){\line(0,1){30}}

\put(40, 845){\line(0,1){40}}
\put(41, 845){\line(0,1){40}}

\put(40, 885){\line(1,0){76}}
\put(40, 886){\line(1,0){76}}

\put(116, 885){\line(0,1){28}}
\put(117, 885){\line(0,1){28}}

\put(116, 912){\line(-1,0){3}}
\put(116, 913){\line(-1,0){3}}

\put(60, 890){$\gamma_n$}
\put(310, 890){$\gamma_n$}

\put(0,800){\circle*{6}}
\put(250,800){\circle*{6}}

\put(98, 839){$z$}
\put(120, 849){${\bf L}_{\theta_p^-}$}
\put(450, 892){${\bf L}_{\theta_p^-}$}


\put(250, 800){\line(0,-1){10}}
\put(251, 801){\line(0,-1){10}}

\put(250, 790){\line(1,0){10}}
\put(251, 791){\line(1,0){10}}

\put(260, 790){\line(0,1){25}}
\put(261, 791){\line(0,1){25}}

\put(260, 815){\line(1,0){30}}
\put(261, 816){\line(1,0){30}}

\put(290, 815){\line(0,1){30}}
\put(291, 816){\line(0,1){30}}

\put(290, 845){\line(0,1){40}}
\put(291, 846){\line(0,1){40}}

\put(290, 885){\line(1,0){80}}
\put(291, 886){\line(1,0){80}}

\put(370, 885){\line(0,-1){8}}
\put(371, 886){\line(0,-1){8}}

\put(370, 876){\line(1,0){30}}
\put(371, 877){\line(1,0){30}}

\put(210, 780){\line(2,1){230}}
\put(310, 730){${\bf D}_{(0,0)} $}

\put(200, 790){\line(2,1){250}}
\put(250, 700){${\bf D}_{(-n,-n)}$}
\put(320, 680){$(b)$}
\put(50, 680){$(a)$}

\put(250, 800){\line(1,-1){50}}

\put(250, 800){\line(-1,1){30}}
\put(220, 770){\line(1,-1){50}}
\put(220, 770){\line(-1,1){30}}
\put(325,950){\line(1,-1){100}}
\put(325,1000){\line(1,-1){130}}
\put(210, 760){\line(2,1){240}}

\put(205,785){\circle*{3}}
\put(210,787){\circle*{3}}
\put(215,789){\circle*{3}}
\put(220,791){\circle*{3}}
\put(225,793){\circle*{3}}
\put(230,795){\circle*{3}}
\put(235,797){\circle*{3}}
\put(240,800){\circle*{3}}
\put(245,803){\circle*{3}}
\put(250,806){\circle*{3}}
\put(255,809){\circle*{3}}
\put(260,812){\circle*{3}}
\put(265,815){\circle*{3}}
\put(270,818){\circle*{3}}
\put(275,820){\circle*{3}}
\put(280,823){\circle*{3}}
\put(285,825){\circle*{3}}
\put(290,828){\circle*{3}}
\put(295,830){\circle*{3}}

\put(300,832){\circle*{3}}
\put(305,835){\circle*{3}}
\put(310,838){\circle*{3}}
\put(315,840){\circle*{3}}
\put(320,842){\circle*{3}}
\put(325,845){\circle*{3}}
\put(330,848){\circle*{3}}
\put(335,850){\circle*{3}}
\put(340,852){\circle*{3}}
\put(345,855){\circle*{3}}
\put(350,858){\circle*{3}}
\put(355,860){\circle*{3}}
\put(360,862){\circle*{3}}
\put(365,865){\circle*{3}}
\put(370,868){\circle*{3}}
\put(375,870){\circle*{3}}
\put(380,872){\circle*{3}}
\put(385,875){\circle*{3}}
\put(390,878){\circle*{3}}
\put(395,880){\circle*{3}}
\put(400,882){\circle*{3}}
\put(405,885){\circle*{3}}
\put(410,888){\circle*{3}}
\put(415,890){\circle*{3}}
\put(420,892){\circle*{3}}
\put(425,895){\circle*{3}}
\put(430,898){\circle*{3}}
\put(220,770){\circle*{3}}
\put(225,772){\circle*{3}}
\put(230,775){\circle*{3}}
\put(235,777){\circle*{3}}
\put(240,780){\circle*{3}}
\put(245,782){\circle*{3}}
\put(250,785){\circle*{3}}

\put(255,787){\circle*{3}}
\put(260,789){\circle*{3}}
\put(265,791){\circle*{3}}
\put(270,793){\circle*{3}}
\put(275,795){\circle*{3}}
\put(280,797){\circle*{3}}
\put(285,799){\circle*{3}}
\put(290,803){\circle*{3}}
\put(295,806){\circle*{3}}
\put(300,809){\circle*{3}}

\put(305,811){\circle*{3}}
\put(310,812){\circle*{3}}
\put(315,815){\circle*{3}}
\put(320,818){\circle*{3}}
\put(325,820){\circle*{3}}
\put(330,822){\circle*{3}}
\put(335,825){\circle*{3}}
\put(340,828){\circle*{3}}
\put(345,831){\circle*{3}}
\put(350,833){\circle*{3}}
\put(355,835){\circle*{3}}
\put(360,838){\circle*{3}}
\put(365,840){\circle*{3}}
\put(370,842){\circle*{3}}
\put(375,845){\circle*{3}}
\put(380,848){\circle*{3}}
\put(385,850){\circle*{3}}
\put(390,852){\circle*{3}}
\put(395,855){\circle*{3}}
\put(400,858){\circle*{3}}
\put(405,860){\circle*{3}}
\put(410,862){\circle*{3}}
\put(415,865){\circle*{3}}
\put(420,868){\circle*{3}}
\put(425,870){\circle*{3}}
\put(430,872){\circle*{3}}
\put(435,875){\circle*{3}}
\put(440,878){\circle*{3}}
\put(445,881){\circle*{3}}

\end{picture}
\end{center}
\caption{{\em The graph in (a) shows how to construct an NE 1-path below ${\bf L}_{\theta_p^-}^- (\epsilon n)$. 
The dotted path is $\gamma$, an NE 1-path between ${\bf L}_{\theta_p^-}$ and ${\bf L}_{\theta_p^-}^+(n^{1/2+\delta})$
from ${\bf D}_{(-n,-n)}$ to ${\bf D}_{(n,n)}$.
The path of small circles  is $\tau_{nu}(\gamma)$, an NE 1-path between $\tau_{nu}({\bf L}_{\theta_p^-})$ and $\tau_{nu}({\bf L}_{\theta_p^-}^+(n^{1/2+\delta}))$
from ${\bf D}_{(0,0)}$ to ${\bf D}_{(2n,2n)}$. The two paths meet at $z$. $\gamma_n$, the boldfaced path,
 is an optimal path
above ${\bf L}_{\theta}^-(\epsilon n)$. It meets $\gamma$ at $v_1$, and it meets $\tau_{nu}(\gamma)$ at $v_2$.
Let us  go along $\gamma_n$ from the origin to $v_1$,  along $\gamma$ from $v_1$ to
$z$,  then  along  $\tau_{nu}(\gamma)$ from $z$ to $v_2$, and finally   along $\gamma_n$ from
$v_2$ to $nu$ to construct an optimal path meeting $L_{\theta}^-(\epsilon n)$. 
The  graph in (b) shows how to construct an
NE 1-path crossing out of ${\bf L}_{\theta_p^-}^+(n^{1/2+\delta})$. The two dotted lines are $\gamma'$ and $\gamma''$
NE 1-paths from ${\bf D}_{(-n,-n)}$ to ${\bf D}_{(2n,2n)}$. An optimal path $\gamma_n$ crosses out of $\gamma'$ at 
$v_1$ and $v_2$, respectively.  Thus, the path along $\gamma_n$ from $v_1$ and $v_2$ has to be an NE 1-path.
}}
\end{figure}

{\bf Proof of Theorem 2.} Recall that ${\bf L}_\theta$ is the line passing through the origin and
$(1, \theta)$. 
The key of the proof for Theorem 2 is
to construct an optimal path such that the distance from some vertices in the path to ${\bf L}_\theta$ 
is $O(n)$. 
We suppose that,  for some $u=(r, \theta)$ for some $\theta_p^- < \theta< \theta_p^+$
$$\xi(u) < 1.\eqno{(3.18)}$$
Without loss of generality, we assume that $r=1$.
Thus, for a large $n$, with a probability larger than $C$, $h_n(u)\leq n^{\xi(u)}$.
In other words, for a small $\epsilon >0$,
all optimal paths from the origin to $nu=(n, \theta)\in {\bf D}_{(n,n)}$ have to stay above the line
${\bf L}_\theta^-(\epsilon n)$ (see Fig. 4). We shall show this is impossible.

Given a graph $G$ on ${\bf Z}^2$, we flip $G$ $180^\circ$ around the line ${\bf D}_{u/2}$, the parallel line in the middle
of ${\bf D}_{(0,0)}$ and ${\bf D}_u$. We denote by $\tau_u(G)$ the new graph.

By Lemma 4, there are an NE $1$-path  $\gamma$ from ${\bf D}_{(-n,-n)}$ to 
${\bf D}_{(n, n)}$ between  ${\bf L}_{\theta_p^-}$ and ${\bf L}_{\theta_p^-}^+(n^{1/2+\delta})$ (see Fig. 4), 
and (by  symmetry)  another NE $1$-path $\tau_{nu}(\gamma)$ between $\tau_{nu}({\bf L}_{\theta_p^-})$ and $\tau_{nu}({\bf L}_{\theta_p^-}^+(n^{1/2+\delta}))$
 with a probability larger than $1-C_1\exp(-C_2 n^{\delta})$ (see Fig. 4).
By our definition, ${\bf L}_{\theta_p^-}^+(n^{1/2+\delta})$ and $\tau_{nu}({\bf L}_{\theta_p^-}^+(2n^{1/2+\delta}))$
meet at ${\bf D}_{(n/2,n/2)}$. We denote by $w$ the intersection point.
 Note that $\theta > \theta_p^-$,
so we take $n$ large 
$$\mbox{{\bf dist}} (w, {\bf L}_\theta) \geq (n -\tan^{-1}(\theta-\theta_p^-)n^{1/2+\delta})\sin(\theta-\theta_p^-)\geq 
n\sin (\theta-\theta_p^-)/2.$$
By this observation, $\gamma$ and $\tau_{nu} (\gamma)$ also meet  at ${\bf D}_{(n/2,n/2)}$. We denote by
$z$ the intersection point. Thus, 
$$\mbox{{\bf dist}}(z, {\bf L}_\theta) \geq \mbox{{\bf dist}} (w, {\bf L}_\theta)\geq n \sin(\theta-\theta_p^-)/2.\eqno{(3.19)}$$

 Now we want to construct an optimal path from the origin to $nu$, but it cannot stay above
${\bf L}_\theta^-(\epsilon n)$ with a large probability. The contradiction tells us that (3.18) cannot hold,
so Theorem 2 follows.
The remaining work is to construct the path.
For a large $n$, by Lemma 4, we assume that $\gamma$ and $\tau_{u/2}(\gamma)$ exist with a probability larger than $1-C/4$.
For a small $\epsilon >0$, we may assume that $\sin(\theta-\theta_p^-)>2\epsilon$.
Therefore,
$\gamma$ meets ${\bf L}_{\theta}^-(\epsilon n)$ at $g$ (see Fig. 4). 
Let us go along from $\gamma$ to $g$,  from $g$ along ${\bf L}_{\theta}^-(\epsilon n)$, back to ${\bf D}_{(-n,-n)}$, and then from ${\bf D}_{(-n,-n)}$ up to $\gamma$. Now we have a triangular shape
enclosing the origin inside the shape.  Note that the shape will stay below ${\bf D}_{(n/2,n/2)}$ as we constructed it.

We select an optimal path $\gamma_n$ from the origin to $nu$. 
As we mentioned, if (3.18) holds, $\gamma_n$ will stay above ${\bf L}_{\theta}^-(\epsilon n)$.
$\gamma_n$ has to cross out of the triangular shape. 
Since $\gamma_n$ is above  ${\bf L}_{\theta}^-(\epsilon n)$,
it crosses out of  either $\gamma$ or ${\bf D}_{(-n,-n)}$. Suppose that
it crosses out of  ${\bf D}_{(-n,-n)}$. Note that $\gamma_n$ has to go back to $nu$ and each edge costs
at least time 1, so $T(\gamma_n) \geq 3n$. We know that 
$$T(\gamma_n)/n \rightarrow \mu n=n \mbox{ a.s.},$$
so if we take  larger $n$,  $\gamma_n$ meets ${\bf D}_{(-n,-n)}$ with a probability less than $C/4$.
Similar, $\gamma_n$ meets ${\bf D}_{(2n,2n)}$ with a probability less than $C/4$.
Now we suppose that $\gamma_n$ meets $\gamma$ at $v_1$ (see Fig. 4). By  symmetry,
$\gamma_n$ also  meets $\tau_{nu}(\gamma)$ at $v_2$ on the existence of $\gamma$ and $\tau_{nu}(\gamma)$.
Note that the path from $v_1$ to $z$, then from $z$ to $v_2$ along $\gamma$ and along
$\tau_{nu}(\gamma)$, respectively, is an NE $1$-path from $v_1$ to $v_2$.
 Therefore,  we go along $\gamma_n$ from the origin to $v_1$,  along $\gamma$ from $v_1$ to
$z$,   along  $\tau_{u/2}(\gamma)$ from $z$ to $v_2$, and finally   along $\gamma_n$ from
$v_2$ to $u$ (see Fig. 4). It is an optimal path.
As our selection,
$$\sin(\theta-\theta_p^-) >2\epsilon,$$
so $z$ has to stay below ${\bf L}_{\theta}^-(\epsilon n)$ with a probability larger than $(1-C/4-C/4-C/4)=1-3C/4$.
It contradicts the assumption that, with a probability larger than $C$, all optimal paths have to stay above ${\bf L}_{\theta}^-(\epsilon n)$. Thus, (3.18) cannot hold, so Theorem 2 follows.\\

{\bf Proof of Theorem 3}. By  symmetry, we only show Theorem 3 when  $u=(r, \theta_p^-)$.
We show all optimal paths staying inside ${\bf S}_{\theta_p^-}(4n^{1/2+\delta})$.
By Lemma  4, there is an NE $1$-path $\gamma'$ from ${\bf D}_{(-n,-n)}^+$ to ${\bf D}_{(2n, 2n)}^+$ between ${\bf L}_{\theta^-_p}^+(n^{1/2+\delta})$ and ${\bf L}_{\theta_p^-}$
with a probability larger than $1-C_1\exp(-C_2 n^{\delta})$ (see Fig. 4).
Similarly, there is another NE $1$-path $\gamma''$ from ${\bf D}_{(-n,-n)}^-$ to ${\bf D}_{(2n, 2n)}^-$ between
 ${\bf L}_{\theta_p^-}$ and ${\bf L}_{\theta^-_p}^-(n^{1/2+\delta})$ with a probability larger than $1-C_1\exp(-C_2 n^{\delta})$ (see Fig. 4).
Then with a probability larger than $1- 2C_1\exp(-C_2 n^{\delta})$, there exist the two 1-paths defined above.
If $h_n(u) \geq n^{1/2+\delta}$, let $\gamma_n$ be an optimal path from the origin to $u$ such that
$v\in \gamma_n$ with a distance from $v$ to $ {\bf L}_{\theta_p^-}$ larger than $n^{1/2+\delta}$. Then either
(a) $\gamma_n$ crosses out of $\gamma'$, or  $\gamma''$, but $\gamma_n$ stays inside
the region between two lines ${\bf D}_{(-n,-n)}$ and  ${\bf D}_{(2n,2n)}$, 
or (b) $\gamma_n$ intersects  the line ${\bf D}_{(-n, -n)}$ or 
${\bf D}_{(2n,2n)}$.
If (a) occurs, suppose that $\gamma_n$ crosses out of $\gamma'$. Let $\gamma_n$ first meet $\gamma'$ at $v_1$. Note that
under (a), $\gamma_n$ will stay inside, between two lines ${\bf D}_{(-n,-n)}$ and  ${\bf D}_{(2n,2n)}$,
so $\gamma_n$ will  meet $\gamma'$ again. Suppose that it meets again at $v_2$  (see Fig. 4).
Note that the path along $\gamma'$ from $v_1$ to $v_2$ is an NE 1-path. 
Note also that $\gamma_n$ is  optimal, so the piece from $v_1$ to
$v_2$ along $\gamma_n$ has to be an NE $1$-path.
 Therefore, we can construct an optimal path along $\gamma'$ to $v_1$ and then along $\gamma_n$ from
$v_1$ to $v_2$, and finally, along $\gamma'$ to $u$. The NE 1-path constructed then will not stay between
${\bf L}_{\theta_p^-}(n^{1/2+\delta})$ and ${\bf L}_{\theta_p^-}$. By Lemma 4, we have
$${\bf P}\left[\mbox{ (a) occurs }\right] \leq C_1 \exp(-C_2 n^{\delta}).\eqno{(3.20)}$$
Similar to our proof for Theorem 2,
$${\bf P}\left[\mbox{ (b) occurs }\right] \leq C_1 \exp(-C_2 n^{\delta}).\eqno{(3.21)}$$
With these observations,
\begin{eqnarray*}
&&{\bf P}\left[h_n(u) > n^{1/2+\delta}\right]\\
&&\leq {\bf P}\left[\mbox{ (a) or (b) occurs}, \exists\,\,\, \gamma'\mbox{ and }\gamma'' \,\,\, \right]+C_1\exp(-C_2 n^{\delta/2})\\
&&\leq 2C_1 \exp(-C_2 n^{\delta}).
\end{eqnarray*}
Therefore, $\xi(u) \leq 0.5$. 

Now we need to show $\xi(u) \geq 0.5$ to complete the proof of Theorem 3. It follows from 
the proof in Lemma 3 of Yukich and Zhang (2006) that, with a positive probability $C$, there is an NE 1-path
$\gamma$ from the origin to $u$ and $\Omega_\infty^{(0,0)}$ occurs. 
Furthermore, $\gamma$ follows from $r_n$, the right hand edge,  at least to ${\bf D}_{(n/2,n/2)}$. As we mentioned in the proof of Lemma 4, 
$r_n$ can be decomposed into an i.i.d. sequence $\{X_i\}$ such that $X_1$ has  an exponential short tail. 
By a standard 
central limit theorem, for all $\delta>0$,
the probability that $r_n$, between ${\bf D}_{(0,0)}$ and ${\bf D}_{(n/2,n/2)}$,
stays inside  strip ${\bf S}_{\theta_p^-}(n^{1/2-\delta})$ goes to zero
as $n\rightarrow \infty$. This shows that with a probability $C$, $h_n(u) \geq n^{1/2-\delta}$ for $\delta >0$.
This shows that $\xi(u) \geq 0.5$,
so Theorem 3 follows.

\section{Proof of Theorem 4.}
By  symmetry, we only need to show Theorem 4 when $u=(r, \theta_p^-)$.
We fix $F(1)=p_0 > \vec{p}_c$. 
Each edge of ${\bf Z}^2$ only takes two values: 1 or a value larger than 1.
Let ${\bf P}_{p_0}(\cdot)$ be the corresponding probability measure, and let ${\bf E}_{p_0}(\cdot)$ be 
the expected value.
Recall that  $(r_{p_0}, \theta_{p_0}^-)$ is  the lower endpoint of the flat segment  (see Fig. 1),
where 
$$r_{p_0}=\sqrt{1/2+ \alpha^2_{p_0}}\mbox{ and }\theta^-_p=\arctan \left( { 1/2-\alpha_p/\sqrt{2} \over  1/2+ \alpha_p/\sqrt{2} }\right ).\eqno{(4.1)}$$

Since we only focus on $(r_{p_0}, \theta_{p_0}^-)$, we replace $\theta^-_{p_0}$ by $\theta_{p_0}$  and $\partial {\bf B}_F$
by $\partial {\bf B}_{p_0}$.
Similarly, for any $p\geq p_0> \vec{p}_c$, 
we define $(r_{p}, \theta_{p})$ as the lower endpoint of the flat segment.
By the Theorem in Zhang (2004), we know that $\alpha_t$ is differentiable at $p$. Also, by (12) in page 1007
of Durrett (1984), we know that 
$$\alpha_p-\alpha_{p_0}\geq 2(p-p_0).\eqno{(4.2)}$$
With (4.1), (4.2), and Zhang's argument, we have the following lemma.\\

{\bf Lemma 5.} {\em If $p_0 > \vec{p}_c$, then $d(\theta_{p})/dp$ exists at $p_0$ and $d(\theta_{p_0})/dp_0<0$}.\\

With Lemma 5, there exists a positive $C=C(p_0)$ such that for $p > p_0$,
$$\theta_{p_0}-\theta_p\leq C(p_0-p).\eqno{(4.3)}$$

For the fixed $p_0$,  we define the following two-variable function for $p\in [p_0, 1]$ and $r >0$.
$$\mu_{p_0}((r, \theta_p))=\lim_{n\rightarrow \infty} {{\bf E}_{p_0}T((0,0), (nr, \theta_p))\over n}.$$
Note that $(r_{p_0}, \theta_{p_0})$ belongs to $\partial {\bf B}_{p_0}$, so
$$\mu_{p_0}((r_{p_0}, \theta_{p_0}))=1.\eqno{(4.4)}$$

For fixed $p_0$ and $\theta_p$, we may consider one variable function
$\mu_{p_0}((r, \theta_p))$ with variable $r$. We show that this function is convex. \\

{\bf Lemma 6.} {\em If $p \geq  p_0 > \vec{p}_c$, for $r_1 >0$ and $r_2 >0$,}
$$\mu_{p_0}((r_1+r_2, \theta_p))\leq \mu_{p_0}((r_1, \theta_p))+ \mu_{p_0}((r_2, \theta_p)).$$

{\bf Proof.} By the subadditive property,
$$T((0,0), (n(r_1+r_2), \theta_p))\leq T((0,0), (nr_1, \theta_p))+T( (nr_1, \theta_p), (n(r_1+r_2), \theta_p)).\eqno{(4.5)}$$
If we take the mean of both sides in (4.5), by  translation invariance, we have
$${\bf E}_{p_0}T((0,0), (n(r_1+r_2), \theta_p))\leq {\bf E}_{p_0}T((0,0), (nr_1, \theta_p))+{\bf E}_{p_0}T((0,0), (nr_2, \theta_p)).
\eqno{(4.6)}$$
Therefore, Lemma 2 follows if we divide by $n$ both sides of (4.5) and let $n\rightarrow \infty$. $\Box$\\

With this convex property, we know that 
$$\mu_{p_0}((r, \theta_p))\mbox{ is continuous in }r\in [0, \infty).\eqno{(4.7)}$$

Now we introduce another lemma, by Yukich and Zhang (2006), to give a more precise estimate.\\

{\bf Lemma 7}  (Yukich and Zhang).  {\em For $\vec{p}\leq p_0 \leq p$, there exists 
a constant $C=C(p_0)$ depending on $p_0$ but not $(p-p_0)$  such that }
$$ \mu_{p_0}((r_p, \theta_p)) \geq 1- {C(p-p_0)^2 \over \log (p-p_0)}.$$

We may select $p\in [p_0, p_0+\delta]$ for some a small $\delta >0$ such that for some $m > 2$
$$ \mu_{p_0}((r_p, \theta_p)) \geq 1+ C(p-p_0)^m.\eqno{(4.8)} $$
With these lemmas, we can show Theorem 4.\\
\begin{figure}\label{F:graph}
\begin{center}
\setlength{\unitlength}{0.0125in}%
\begin{picture}(200,230)(67,800)
\thicklines
\put(0, 800){\vector(1,0){300}}
\put(0, 800){\vector(0,1){300}}
\put(75,950){\circle*{4}}
\put(150,875){\circle*{4}}
\put(0, 800){\line(1,2){100}}
\put(0, 800){\line(2,1){200}}
\put(75,950){\line(1,-1){100}}
\put(200,900){${{\theta_{p_0}}}$}
\put(145,890){${{r_{p_0}}}$}
\put(0, 800){\line(3,1){200}}
\put(200,860){${{\theta_{p}}}$}
\put(167,870){${{d^+_\theta}}$}

\put(120,910){${{\bf S}^+_{\theta_0}}$}
\put(165,825){\circle*{3}}

\put(153,865){\circle*{3}}
\put(156,853){\circle*{4}}
\put(159,845){\circle*{3}}
\put(162,835){\circle*{3}}
\put(150,875){\circle*{3}}
\put(140,885){\circle*{3}}
\put(130,895){\circle*{3}}
\put(120,905){\circle*{3}}
\put(110,915){\circle*{3}}
\put(100,925){\circle*{3}}
\put(142,835){$\bar{r}_p$}
\put(167,856){\circle*{4}}
\put(180,847){$r_p$}

\put(167,825){$\partial {\bf B}_{p_0}$}

\end{picture}
\end{center}
\caption{{\em The graph shows that there is left $l$-curvature at $(r_{\theta_0}, \theta_0)$.
The dotted line is $\partial {\bf B}_{p_0}$.}}
\end{figure}
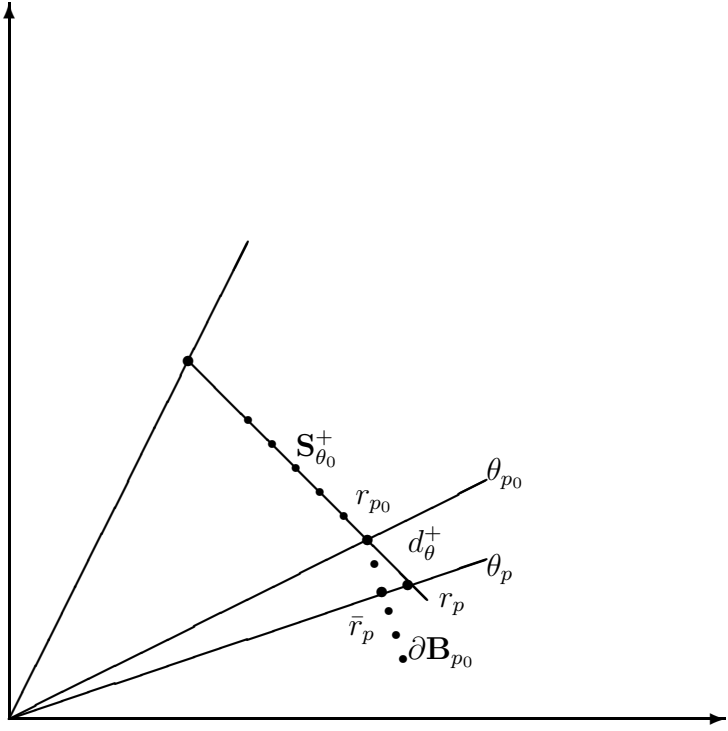

{\bf Proof of Theorem 4.} Note that $\mu((0, \theta_p))=0$, so by (4.7), (4.8), and  the intermediate value theorem, there exists $\bar{r}_p$,
for $p > p_0$
$$\mu_{p_0}(\bar{r}_p, \theta_p)=1\eqno{(4.9)}$$ 
for $\bar{r}_p <  r_p$.
By the shape theorem, 
$(\bar{r}_p, \theta_p)\in \partial {\bf B}_{p_0}$.
We choose a number $q\geq 0$ such that $\bar {r}_p= r_{p}-q$.
By using the convex  property in Lemma 6 and (4.8), we have for $m >2$,
$$
\mu_{p_0}((q, \theta_p))+\mu_{p_0}( (\bar{r}_p,\theta_p))\geq \mu_{p_0}((r_{p},\theta_p))\geq 1+ 
(p-p_0)^m .\eqno{(4.10)}$$
Thus, we have by (4.9)
$$
\mu_{p_0}((q,\theta_p))\geq (p-p_0)^m. \eqno{(4.11)}
$$
Note that there is a path that is at most $2n+1$ edges from $(0,0)$ to
$(nq, \theta_p)$ (if it is in ${\bf Z^2}$) or to the integer valued coordinate next to $(nq,\theta_p)$, so 
$$T((0,0), (nq, \theta_p)) \leq 2n q+2 .\eqno{(4.12)}$$
By (4.12),
$$\mu_{p_0}((q, \theta_p))\leq 2q. \eqno{(4.13)}$$
With this observation together with (4.11),
 $$r_p- \bar{r}_p =q \geq \mu_{p_0}((q, \theta_p))/2 \geq (p-p_0)^m/2.\eqno{(4.14)}$$
By (4.3), there exists $C=C(p_0)$ such that
$$r_p-\bar{r}_p\geq C(\theta_{p_0}-\theta_p)^m.\eqno{(4.15)}$$
By our definition, we have
$$(r_p, \theta_p)=\left(\sqrt{1/2+\alpha^2_p},\theta_p\right)\in \{(x,y): x+y=1\}.\eqno{(4.16)}$$
On the other hand, ${\bf S}^+_{\theta_0}$, defined in section 1, is the line $ \{(x,y): x+y=1\}$.
When $\theta_p$ increases to $\theta_{p_0}$, the boundary curve along  $\partial {\bf B}_{p_0}$ 
(see Fig. 5) goes to $(r_{p_0}, \theta_{p_0})$, but  below  $\{x+y=1\}$ (see Fig. 5). 
By  (4.15) and a simple computation, we have for some $C$ and $m > 2$,
$$\mbox{{\bf dist}}\left((\bar{r}_{p}, \theta\right), {\bf S}^+_{\theta_0})\geq  
 (r_p- \bar{r}_p)\geq C(\theta_{p_0}-\theta_p)^m\geq C \left[(d^+_\theta)\sin(\pi/4)\right]^m.\eqno{(4.17)}$$
By (4.17) and the definition of the right $l$-curvature, we have
$$\kappa^+((r_{\theta_0}, \theta_0))\geq 0.5.\eqno{(4.18)}$$
Therefore, Theorem 4 follows from (4.18).

\begin{center}
{\bf References}
\end{center}
Alexander, K. (1997). Approximation of subadditive functions and convergence rates in limiting-shape results. {\em Ann.  Probab.} {\bf 25}
30--55.\\
Cox, T. and Durrett, R. (1981) Some limit theorems for percolation processes with necessary and sufficient conditions. {\em Ann. Probab.} {\bf 9} 583--603.\\
Durrett, R. (1984). Oriented percolation in two
dimensions. {\em Ann. Probab.} {\bf 12} 999--1040.\\
Durrett, R. and  Liggett, T. (1981).
The shape of the limit set in Richardson's growth model. {\em Ann.
Probab.} {\bf 9} 186--193.\\
Grimmett, G. (1999). {\em Percolation.} Springer, Berlin.\\
 Hammersley, J. M. and  Welsh, D. J. A. (1965).
First-passage percolation, subadditive processes,
stochastic networks and generalized renewal theory.
In {\em Bernoulli, Bayes, Laplace Anniversary Volume} 
(J. Neyman   and L. LeCam,  eds.) 61--110. Springer, Berlin.\\
Johansson, K. (2000). Transversal fluctuations for increasing subsequences on the plane.
{\em PTRF} {\bf 116} 445--456.\\
Kesten, H. (1986). Aspects of first-passage percolation. {\em Lecture Notes in
Math.} {\bf 1180} 125--264. Springer, Berlin.\\
Kesten, H. (1993). On the speed of convergence in first passage percolation. {\em Ann. Appl. Probab.} {\bf 3} 296--338.\\
Kingman, J.F.C. (1973). Subadditive ergodic theory. {\em Ann. Probab.} {\bf 1} 883--909.\\
Krug, J. and Spohn, H. (1992). Kinetic roughening of growing surfaces. In {\em Solids Far from
Equilibrium: Growth, Morphology, Defects} (C. Godreche, ed.) 497--582. Cambridge Univ. Press, Cambridge.\\
Kuczek, T. (1989). The central limit theorem for the
right edge of supercritical oriented percolation. {\em  Ann. Probab.} {\bf
17} 1322--1332.\\
Howard, D. (2000). Models of first-passage percolation. In {\em Probability on discrete structures} (Kesten, H. eds)
125--174. Springer, Berlin.\\
Marchand, R. (2002). Strict inequalities for the time
constant in first passage percolation. {\em Ann. Appl. Probab.} {\bf 12}
1001--1038.\\
Newman, C. and Piza, M. (1995). Divergence of shape fluctuations in two dimensions. {\em Ann. Probab.} {\bf 23} 977--1005. \\
Smythe, R. T. and Wierman, J. C. (1978).
First passage percolation on the square lattice.
{\em Lecture Notes in Math.} {\bf 671}. Springer, Berlin.\\
Yukich, J. and Zhang, Y. (2006). Singularity points for first passage percolation. {\em Ann. Probab.} {\bf 34} 577--592. \\
Zhang, Y. (2004).
On the infinite differentiability of the right edge in the
supercritical oriented percolation, {\em Stoch. Proc. Appl.} {\bf 114} 279--286.\\
Zhang, Y (2006). Shape fluctuations are different in  different directions. (To appear {\em Ann. Probab}).

\noindent
Yu Zhang\\
Department of Mathematics\\
University of Colorado\\
Colorado Springs, CO 80933\\
email: yzhang3@uccs.edu\\

\end{document}